\documentclass[11pt,twoside]{amsart}

\usepackage{enumerate}
\usepackage[hyperindex=true,plainpages=false,colorlinks=false,pdfpagelabels]{hyperref}\DeclareMathOperator{\End}{End} \DeclareMathOperator{\Hom}{Hom}
 \DeclareMathOperator{\SL}{SL}
 
\DeclareMathOperator{\Span}{Span} \DeclareMathOperator{\Id}{Id}
 \DeclareMathOperator{\diam}{diam}
\DeclareMathOperator{\dist}{dist} 
\DeclareMathOperator{\Harm}{Harm}
\DeclareMathOperator{\Spec}{Spec}
\DeclareMathOperator{\Specs}{\widetilde\Spec}
\DeclareMathOperator{\SpeC}{\mathcal{S}}
\DeclareMathOperator{\SpeCs}{\,\widetilde{\! \SpeC}}

\DeclareMathOperator{\Ind}{Index}

\renewcommand{\ker}{\operatorname{ker}}
\newcommand{\im}{\operatorname{im}}

\newcommand{\coker}{\operatorname{coker}}
\newcommand{\Index}{\operatorname{index}}

\DeclareMathOperator{\dbar}{\bar\partial}

\DeclareMathOperator{\ord}{ord}

\newcommand{\R}{\mathbb{R}}
\newcommand{\C}{\mathbb{C}}

\newcommand{\Z}{\mathbb{Z}}
\renewcommand{\P}{\mathbb{P}}
\renewcommand{\H}{\mathbb{H}}  

\newcommand{\CP}{\mathbb{CP}}
\newcommand{\HP}{\mathbb{HP}}
\renewcommand{\S}{S}

\newcommand{\eprint}[1]{e--print: \href{http://#1}{\nolinkurl{#1}}}

\theoremstyle{plain}
\newtheorem{The}{Theorem}[section]
\newtheorem{Pro}[The]{Proposition}
\newtheorem{Lem}[The]{Lemma}
\newtheorem{Cor}[The]{Corollary}
\newtheorem*{Cor*}{Corollary}

\theoremstyle{definition}
\newtheorem{Def}[The]{Definition}

\theoremstyle{remark} 
\newtheorem{Rem}[The]{Remark}

\newtheorem*{Rem*}{Remark}

\numberwithin{equation}{section}
\begin{document}

\title[The spectral curve of a quaternionic holomorphic line bundle over
$T^2$]{The spectral curve of a quaternionic holomorphic line bundle over a
  2--torus}

\date{\today}

\author{Christoph Bohle}
\author{Franz Pedit}
\author{Ulrich Pinkall}

\address{Christoph Bohle\\
  Institut f\"ur Mathematik\\
  Technische Universit{\"a}t Berlin\\
  Stra{\ss}e des 17.\ Juni 136\\
  10623 Berlin\\
  Germany}

\address{Franz Pedit\\
  Mathematisches Institut der  Universit{\"a}t T\"ubingen\\
  Auf der Morgenstelle 10\\
  72076 T\"ubingen\\
  and\\
  Department of Mathematics \\
  University of Massachusetts\\
  Amherst, MA 01003, USA }

\address{Ulrich Pinkall\\
  Institut f\"ur Mathematik\\
  Technische Universit{\"a}t Berlin\\
  Stra{\ss}e des 17.\ Juni 136\\
  10623 Berlin\\
  Germany}

\email{bohle@math.tu-berlin.de,
  pedit@mathematik.uni-tuebingen.de,\phantom{grrrrrrrrrrrrrrrr}\linebreak
  pinkall@math.tu-berlin.de}

\thanks{All authors supported by DFG SPP 1154 ``Global Differential
  Geometry''. Second author additionally supported by Alexander von Humboldt foundation}

\maketitle

\section{Introduction}\label{sec:intro}
Over the last 20 years algebraically completely integrable systems have been
studied in a variety of contexts.  On the one hand their theory is
interesting from a purely algebraic geometric point of view, on the other hand
a number of problems arising in mathematical physics and global differential
geometry can be understood in the framework of those integrable systems.  This
situation has led to a rich cross--fertilization of algebraic geometry, global
differential geometry and mathematical physics.

The phase space of an algebro--geometric integrable system consists of moduli
of algebraic curves together with their Jacobians, the Lagrangian tori on
which the motion of the system linearizes in a direction osculating the Abel
image (at some marked point) of the curve.  In classical terminology the
algebraic curve, usually referred to as {\em the spectral curve}, encodes the
action variables whereas the Jacobian of the curve encodes the angle variables
of the system. Since the Jacobian of a curve acts on the Picard variety of fixed
degree line bundles, the linear flow in the Jacobian gives rise, via the
Kodaira embedding, to a flow of algebraic curves in a projective space.  If
this projective space is $\P^1$ the correct choice of moduli of curves yields
as flows harmonic maps from $\R^2$ to $\P^1=S^2$. In case those flows are
periodic in the Jacobian, one obtains harmonic cylinders or tori in $S^2$.
The action of the Jacobian on harmonic maps is a geometric manifestation of
the sinh--Gordon hierarchy in mathematical physics, where one perhaps is more
interested in the solutions of the field equation, here the elliptic
sinh--Gordon equation, rather than the harmonic maps described by these
solutions.  It is a remarkable fact \cite{PS89,Hi} that, due to the
ellipticity of the harmonic map equation, all harmonic tori arise in this
way. It is well known that harmonic maps into $S^2$ are the unit normal maps
of constant mean curvature surfaces in $3$--space.  Thus, the classical
problem of finding all constant mean curvature tori can be rephrased by
studying a particular algebro-geometric integrable system \cite{PS89,Bob}.

This picture pertains in many other instances, including finite gap solutions
of the KdV hierarchy, whose ``fields'' can be interpreted as the Schwarzian
derivatives of curves in $\P^1$, elliptic Toda field equations for the linear
groups, which arise in the study of minimal tori in spheres and projective
spaces, and Willmore tori in $3$ and $4$--space. In each of these cases the
equations can be rephrased as an algebraically completely integrable system.
One of the implications of such a description is the explicit computability of
solutions, since linear flows on Jacobians are parametrized by theta functions
of the underlying algebraic curve.  Moreover, the energy functional ---
Dirichlet energy, area, Willmore energy etc. --- of the corresponding
variational problem appears as a residue of a certain meromorphic form on the
spectral curve, which makes the functional amenable to algebro--geometric
techniques.

In all of the above cases the spectral curves are intimately linked to
holonomy representations of holomorphic families of flat connections. For
example, the harmonic map equation of a Riemann surface into the 2-sphere can
be written as a zero curvature condition on a $\C_*$--family of
$\SL(2,\C)$--connections \cite{Hi,FLPP01}. If the underlying surface is a
2--torus, the eigenvalues of the holonomy (based at some point on the torus)
give a hyper--elliptic curve.  The corresponding eigenlines define a line
bundle over this curve which moves linearly in the Picard of the curve as the
base point moves over the torus.  Perhaps there has been a sentiment that this
setup is prototypical for algebraically completely integrable systems which
arise in the context of differential geometry.  We are learning
now~\cite{BLPP} that there is a much less confined setting for which the above
described techniques are applicable: the geometric classes of surfaces
described by zero curvature equations are special invariant subspaces of a
more general phase space related to the Davey--Stewartson hierarchy
\cite{Ko96,Ta97,Ko00,BPP01,McI}.

This more general setting has to do with the notion of a ``spectral variety''
\cite{DKN,N,Kr,Ku93,FKT03} for a differential operator.  For example, if we
want to study the spectrum of the Schr\"odinger operator on a periodic
structure (crystal), we need to find a solution (the wave function) for a
given energy which is quasi-periodic, that is, gains a phase factor over the
crystal, since the physical state only depends on the complex line spanned by
the wave function.  If the crystal is a 2--D lattice, we obtain a spectral
variety given by the energies and the possible phases of the wave functions,
which in this case would be an analytic surface. Another example occurs
in the computations of tau and correlation functions of massive conformal
field theories over a Riemann surface $M$. In this case one is interested in
solutions with monodromy of the Dirac operator with mass
\[
D=\begin{pmatrix} \dbar &-m\\
m & \partial
\end{pmatrix}
\]
that is, solutions $\psi$ which satisfy
\[
D\psi=0\quad\text{and}\quad \gamma^*\psi =\psi h_{\gamma}
\]
for a representation $h\colon\pi_1(M)\to\C_*$ of the fundamental group acting
by deck transformations.  The spectral variety, parametrizing the possible
monodromies, generally is an analytic set.

This last example arises naturally \cite{FLPP01} in the study of quaternionic
holomorphic line bundles $W$ over a Riemann surface $M$. Such lines bundles
carry a complex structure $J\in\Gamma(\End(W))$ compatible with the
quaternionic structure and the holomorphic structure is described by a
quaternionic -- generally not complex -- linear first order operator
\[
D=\dbar +Q\colon \Gamma(W)\to \Gamma(\bar{K}W)\,.
\]
Here $Q\in\Gamma(\bar{K}\End_{-}(W))$ is the complex anti--linear part of $D$
and the complex linear part $\dbar$ is a complex holomorphic structure on $W$.
The energy of the holomorphic line bundle $W$ is the $L^2$--norm
\[
\mathcal{W}(W,D)=2\int_{M}<Q\wedge * Q>
\]
of $Q$ which is zero for complex holomorphic structures $D=\dbar$. There are
two important geometric applications depending on the dimension $h^0(W)$ of
the space of holomorphic sections $H^0(W)$.  If $h^0(W)=1$ and the spanning
holomorphic section $\psi\in H^0(W)$ has no zeros, the complex structure $J$
can be regarded as a smooth map from the Riemann surface $M$ to $S^2$ whose
Dirichlet energy is $\mathcal{W}$. In the case $h^0(W)=2$ the ratio of two
independent holomorphic sections defines a (branched) conformal immersion from
the Riemann surface $M$ to $S^4$ whose Willmore energy, the average square
mean curvature $\int_{M} H^2$, is given by~$\mathcal{W}$. In both cases all
the respective maps are described by suitably induced quaternionic holomorphic
line bundles \cite{FLPP01,BFLPP02}.

A central geometric feature of both theories is the existence of Darboux
transformations \cite{BLPP} which is a first manifestation of complete
integrability.  These Darboux transforms correspond to 
holomorphic sections with monodromy, that is, sections $\psi\in
H^0(\tilde{W})$ of the pullback bundle $\tilde W$ of $W$ to the universal
cover of $M$, which satisfy
\[ \gamma^*\psi=\psi h_{\gamma}\] for a representation $h\colon \pi_1(M)\to
\H_{*}$.

From here on we only discuss the case when the underlying Riemann surface $M$
is a $2$--torus $T^2=\R^2/\Gamma$. Then the fundamental group is abelian and it
suffices to consider holomorphic sections with complex monodromies which are
all of the form $h=\exp(\int\omega)$ for a harmonic form $\omega\in
\Harm(T^2,\C)$.  Interpreting $\exp(\int\omega)$ as a (non--periodic) gauge on
$\R^2$, a holomorphic section $\psi$ with monodromy $h=\exp(\int\omega)$ gives
rise to the section $\psi\exp(-\int\omega)\in\Gamma(W)$ without monodromy on
the torus which lies in the kernel of the periodic operator
\[
D_{\omega}=\exp(-\int\omega)\circ D\circ \exp(\int\omega)
\colon \Gamma(W)\to \Gamma(\bar{K}W)
\,.
\]
Therefore, the Darboux transforms are described by the harmonic forms
$\omega\in \Harm(T^2,\C)$ for which $D_{\omega}$ has a non--trivial kernel,
and we call this set the {\em logarithmic spectrum} $\Specs(W,D)$ of the
quaternionic holomorphic line bundle $W$.  Note that $\Specs(W,D)$ is
invariant under translations by the dual lattice $\Gamma^*$ of integer period
harmonic forms and its quotient under this lattice, the {\em spectrum}
\[
\Spec(W,D)=\Specs(W,D)/\Gamma^*\subset \Hom(\Gamma,\C_*)\,,
\] 
is the set of possible monodromies of the holomorphic structure $D$.  The
spectrum carries a real structure $\rho$ induced by complex conjugation:
 if a section $\psi$ has monodromy~$h$ the section $\psi j$ has
monodromy $\bar{h}$.

In the description of surface geometry via solutions to a Dirac operator with
potential the spectrum already appeared in the papers of Taimanov~\cite{Ta98},
and Grinevich and Schmidt~\cite{GS98}, and its relevance to the Willmore
problem can be seen in \cite{S02,Ta06} and \cite{BLPP, B}.



The present paper  analyzes  
the structure of the spectrum $\Spec(W,D)$ of a quaternionic line bundle $W$
with holomorphic structure $D$ of degree zero over a 2--torus.  Due to
ellipticity the family $D_\omega$, parametrized over harmonic 1--forms
$\Harm(T^2,\C)$, is a holomorphic family of Fredholm operators.  The minimal
kernel dimension of such a family is generic and attained on the complement of
an analytic subset. In Sections~\ref{sec:foundation} and~\ref{sec:analysis} we
show that for a degree zero bundle over a 2--torus these operators have
$\Index(D_\omega)=0$, their generic kernel dimension is zero, and the spectrum
$\Spec(W,D)\subset \Hom(\Gamma,\C_*)$ is a 1--dimensional analytic set. The
spectrum can therefore be normalized $h\colon \Sigma\rightarrow \Spec(W,D)$ to
a Riemann surface $\Sigma$, the {\em spectral curve}.  Moreover, for generic
$\omega\in \widetilde\Spec(W,D)$ the kernel of $D_\omega$ is 1--dimensional
and therefore $\ker(D_\omega)$ gives rise to a holomorphic line bundle
$\mathcal{L}$, the {\em kernel bundle}, over the spectral curve $\Sigma$.  The
fiber of $\mathcal{L}$ over a generic point $\sigma \in \Sigma$ is the space
of holomorphic sections of $W$ with monodromy $h^\sigma\in\Spec(W,D)$.  The
real structure $\rho$ on the spectral curve $\Sigma$ induced by $\rho(h)=\bar
h$ is covered by multiplication by $j$ on the kernel bundle $\mathcal{L}$ and
therefore has no fixed points.

Physical intuition suggests that for large monodromies the spectrum should be
asymptotic to the {\em vacuum} spectrum $\Spec(W,\dbar)$. That this is indeed
the case we show in Section~\ref{sec:asymptotic_geometry}. Since the vacuum is
described by complex holomorphic sections with monodromy it is a translate of
the curve
\[
\exp(H^0(K))\, \cup\, \exp(\overline{H^0(K)})\subset \Hom(\Gamma,\C_*)
\]
with double points along the lattice of real representations.  We show that
outside a sufficiently large compact subset of $ \Hom(\Gamma,\C_*)$ the
spectrum is a graph over the vacuum, at least away from the double
points. Near a double point the spectrum can have a handle.  Depending on
whether an infinite or finite number of handles appear, the spectral curve has
infinite genus, is connected and has one end, or its genus is finite, it has
two ends and at most two components.

For a generic holomorphic line bundle $W$ the spectral curve $\Sigma$ will
have infinite genus and algebro--geometric techniques cannot be applied. This
motivates us to study the case of finite spectral genus in
Section~\ref{sec:finite_genus} in more detail.  The end behavior of the
spectrum then implies that outside a sufficiently large compact set in
$\Hom(\Gamma,\C_*)$ none of the double points of the vacuum get resolved into
handles. Therefore, we can compactify the spectral curve~$\Sigma$ by adding
two points $o$ and $\infty$ at infinity which are interchanged by the real
structure~$\rho$.  Because the kernel bundle is asymptotic to the kernel
bundle of the vacuum, the generating section $\psi^\sigma\in H^0(\tilde{W})$
with monodromy $h^\sigma$ of $\mathcal{L}_\sigma$ has no zeros for $\sigma$ in
a neighborhood of $o$ or $\infty$.  We show that the complex structures
$S^\sigma$ defined via $S^\sigma\psi^\sigma=\psi^\sigma i$ are a holomorphic
family limiting to $\pm J$ when $\sigma$ tends to $o$ and $\infty$, where $J$
is the complex structure of our quaternionic holomorphic line bundle
$W$. Thus, $S$ extends to a $T^2$-family of algebraic functions
$\overline{\Sigma}\rightarrow \CP^1$ on the compactification
$\overline{\Sigma}=\Sigma\cup\{o,\infty\}$. Pulling back the tautological
bundle over $\CP^1$ the family $S$ gives rise to a linear $T^2$--flow of line
bundles in the Picard group of $\overline{\Sigma}$.  Finally we give a formula
for the Willmore energy~$\mathcal{W}$ in terms of the residue of the
logarithmic derivative of $h$ and the conformal structure of $T^2$.  This
formula allows various interpretations of the Willmore energy including an
interpretation as the convergence speed of the spectrum to the vacuum
spectrum.

During the preparation of this paper the authors profited from conversations
with Martin Schmidt and Iskander Taimanov.

\section{The Spectrum of a Quaternionic Holomorphic Line
  Bundle}\label{sec:foundation}
  
In this section we summarize the basic notions of quaternionic holomorphic
geometry ~\cite{FLPP01} in as much as they are relevant for the purposes of
this paper.  We also recall the basic definitions and properties of the
spectrum of a holomorphic line bundle which, from a surface geometric point of
view, can also be found in ~\cite{BLPP}.

\subsection{Preliminaries} Let $W$ be a quaternionic (right) vector bundle
over a Riemann surface $M$.  A complex structure $J$ on $W$ is a section $J\in
\Gamma(\End(W))$ with $J^2=-\Id$, or, equivalently, a decomposition $W = W_+
\oplus W_-$ into real subbundles which are invariant under multiplication by
the quaternion $i$ and interchanged by multiplication with the quaternion $j$:
the $\pm i$--eigenbundle of~$J$ is $W_{\pm}$.  Note that $W_{+}$ and $W_{-}$
are isomorphic via multiplication by $j$ as vector bundles with complex
structures $J_{|_{W_{\pm} }}$.  The degree of the quaternionic bundle $W$ with
complex structure $J$ is defined as the degree $\deg W := \deg W_+$ of the
underlying complex vector bundle $W_+$.

A {\em quaternionic holomorphic} structure on a quaternionic vector bundle $W$
equipped with a complex structure $J$ is a quaternionic linear operator
\begin{equation*}
D\colon \Gamma(W)\to\Gamma(\bar{K}W)
\end{equation*}
that satisfies the Leibniz rule $D(\psi\lambda) = (D\psi)\lambda + (\psi
d\lambda)''$ for all $\psi\in\Gamma(W)$ and $\lambda\colon M \to \H$ where
$(\psi d\lambda)''=\frac12(\psi d\lambda+ J*\psi d\lambda)$ denotes the $(0,1)$--part of the $W$--valued $1$--form $\psi d\lambda$ with respect to the complex structures $J$ on $W$ and $*$ on $TM^{*}$.  Decomposing
the operator $D$ into $J$ commuting and anti--commuting parts gives $D = \dbar
+ Q$ where $\dbar=\dbar\oplus\dbar$ is the double of a complex holomorphic
structure on $W_+$ and $Q\in\Gamma(\bar{K}\End_{-}(W))$ is a $(0,1)$--form 
with values in the complex anti--linear endomorphisms of $W$ (with complex structure post--composition by $J$), called the \emph{Hopf field}.

The space of holomorphic sections of $W$ is denoted by $H^{0}(W)= \ker (D)$.
Because $D$ is an elliptic operator its kernel $H^0(W)$ is finite dimensional
if the underlying surface is compact.  The $L^2$--norm
\begin{equation*}
\label{eq:Willmore_energy_holbundle}
\mathcal{W}(W)= \mathcal{W}(W, D) = 2\int_M <Q\wedge *Q>
\end{equation*}
of the Hopf field $Q$ is called the {\em Willmore energy} of the holomorphic
bundle $W$ where $<\,,\,>$ denotes the real trace pairing on $\End(W)$. The special
case $Q=0$, for which $\mathcal{W}(W)=0$, describes (doubles of) complex
holomorphic bundles $W\cong W_+\oplus W_+$.

\subsection{The quaternionic spectrum of quaternionic holomorphic line
  bundles}
A holomorphic structure $D$ on $W$ induces a quaternionic holomorphic structure on the pullback bundle
$\tilde W=\pi^* W$ by the universal covering $\pi\colon \tilde M\rightarrow M$.
The operator $D$ on $\tilde{W}$ is periodic with respect to the group $\Gamma$ of deck
transformations of $\pi\colon \tilde M\rightarrow M$. The space $H^0(W)$ of
holomorphic sections of $W$ are the periodic, that is, $\Gamma$--invariant,
sections of $\tilde W$ solving $D\psi=0$.  In the following we also need to consider 
solutions with monodromy of $D\psi=0$. Such solutions
are the holomorphic sections of $\tilde
W$ that satisfy
\begin{equation}
  \label{eq:monodromy}
  \gamma^*\psi=\psi h_\gamma \qquad \textrm{ for all } \gamma\in \Gamma,
\end{equation}
where $h\in \Hom(\Gamma,\H_*)$ is a representation of $\Gamma$. A holomorphic
section $\psi\in H^0(\tilde W)$ satisfying \eqref{eq:monodromy} for some $h\in
\Hom(\Gamma,\H_*)$ is called a \emph{holomorphic section with monodromy} $h$
of $W$, and we denote the space of all such sections by $H^0_h(\tilde W)$.  By
the quaternionic Pl\"ucker formula with monodromy (see appendix
to~\cite{BLPP}) $H^0_h(\tilde W)$ is a finite dimensional real vector
space. Multiplying $\psi\in H^0_h(\tilde W)$ by some $\lambda\in \H_*$ yields
the section $\psi\lambda$ with monodromy $\lambda^{-1} h \lambda$.  Unless $h$
is a real representation $H^0_h(\tilde W)$ is not a quaternionic vector space.

\begin{Def}\label{def:spec}
  Let $W$ be a quaternionic line bundle with holomorphic structure $D$ over a
  Riemann surface $M$.  The \emph{quaternionic spectrum} of $W$ is the
  subspace
  \begin{equation*}
    \Spec_\H(W,D) \subset \Hom(\Gamma,\H_*)/\H_*
  \end{equation*}
  of conjugacy classes of monodromy representations occurring for
  holomorphic sections of $\tilde{W}$. In other words, $h$ represents a point
  in $\Spec_\H(W,D)$ if and only if there is a non--trivial holomorphic
  section $\psi\in H^0(\tilde{W})$ with monodromy $h$, that is, a solution of
  \[  
  D\psi = 0\, \quad\text{satisfying }\quad \gamma^*\psi=\psi h_\gamma \quad
  \textrm{ for all } \gamma\in \Gamma.
  \]
\end{Def}

Our principal motivation for studying the quaternionic spectrum is the
observation \cite{BLPP} that Darboux transforms of a conformal
immersion $f\colon M\rightarrow \S^4$ correspond to nowhere vanishing
holomorphic sections with monodromy of a certain quaternionic holomorphic line bundle
induced by $f$.  From this point of view the quaternionic
spectrum arises as a parameter space for the space of Darboux
transforms.  The quaternionic holomorphic line bundle induced by  $f\colon
M\rightarrow \S^4$ is best described  when viewing $S^4=\HP^1$ as the quaternionic projective line.
The immersion $f$ is the pull--back $L$ of the tautological line bundle over $\HP^1$ 
and as such a subbundle $L\subset \underline{\H^2}$ of the trivial quaternionic rank~2
bundle. The induced quaternionic
holomorphic line bundle is the quotient bundle $W=V/L$ equipped with the
unique holomorphic structure for which constant sections of $V$
project to holomorphic sections of $V/L$.

The idea of defining spectra of conformal immersions first appears, for tori
in $\R^3$, in the work of Taimanov \cite{Ta98}, and Grinevich and Schmidt
\cite{GS98}.  Their definition leads to the same notion of spectrum although
it is based on a different quaternionic holomorphic line bundle associated
to a conformal immersion into Euclidean 4--space $\R^4$ via the
Weierstrass representation \cite{PP}.

\subsection{The spectrum of a quaternionic holomorphic line bundle over
  a 2--torus}
From here on we study the geometry of $\Spec_\H(W,D)$ in the case that the
underlying surface is a torus $T^2=\R^2/\Gamma$. Due to the abelian
fundamental group $\pi_1(T^2)=\Gamma$ every representation
$h\in\Hom(\Gamma,\H_*)$ of the group of deck transformations can be conjugated
into a complex representation in $\Hom(\Gamma,\C_*)$.  The complex
representation $h\in\Hom(\Gamma,\C_*)$ in a conjugacy class in
$\Hom(\Gamma,\H_*)$ is uniquely determined up to complex conjugation $h\mapsto
\bar h$ (which corresponds to conjugation by the quaternion $j$). In
particular, away from the real representations, the map
\[
\Hom(\Gamma,\C_*)\to\Hom(\Gamma,\H_*)/\H_* 
\]
is $2:1$. The lift of the quaternionic spectrum $\Spec_\H(W,D)$ of $W$
under this map gives rise to the spectrum of the quaternionic holomorphic
line bundle $W$.
\begin{Def}
  Let $W$ be a quaternionic holomorphic line bundle over the torus with
  holomorphic structure $D$. Its \emph{spectrum} is the subspace
\[
 \Spec(W,D)\subset\Hom(\Gamma,\C_*)
\]
of complex monodromies occurring for non--trivial holomorphic sections
of~$\tilde W$.
\end{Def}
By construction, the spectrum is invariant under complex conjugation
$\rho(h)=\bar{h}$ and
\begin{equation}
\label{eq:mod-tau}
\Spec_\H(W,D)=\Spec(W,D)/\rho.
\end{equation}

The study of $\Spec(W,D)$ is greatly simplified by the fact that
$G=\Hom(\Gamma,\C_*)$ is an abelian Lie group with Lie algebra
$\mathfrak{g}=\Hom(\Gamma,\C)$ whose exponential map $\exp\colon
\mathfrak{g}\rightarrow G$ is induced by the exponential function
$\C\rightarrow \C_*$.  The Lie algebra $\Hom(\Gamma,\C)$ is isomorphic to the
space $\Harm(T^2,\C)$ of harmonic 1--forms: a harmonic 1--form $\omega$ gives
rise to the period homomorphism $\gamma\in \Gamma \mapsto \int_\gamma \omega$. The
image under the exponential map of such a homomorphism is the multiplier
$\gamma\mapsto h_\gamma= e^{\int_\gamma\omega}$.  The kernel of the group
homomorphism $\exp$ is the lattice of harmonic forms $\Gamma^* =
\Harm(T^2,2\pi i \Z)$ with integer periods. The exponential function thus
induces an isomorphism
\begin{equation}
\label{eq:isomorphism}
\Harm(T^2, \C)/\Gamma^* \cong \Hom(\Gamma,\C_*)\,.
\end{equation}

Rather than $\Spec(W,D)$ we study its logarithmic image which is the
$\Gamma^*$--invariant subset $\Specs(W,D) \subset \Harm(T^2,\C)$, the {\em
  logarithmic spectrum}, with the property that
\begin{equation}
\label{eq:spect_tilde}
\Spec(W,D) \cong\Specs(W,D)/\Gamma^*.
\end{equation}
Considering $\Specs(W,D)$ has the advantage that it is described as the
locus of those $\omega\in \Harm(T^2,\C)$ for which the operator 
\begin{equation}
\label{eq:D_omega_def}
D_{\omega}=e^{-\int \omega} \circ D\circ e^{\int\omega}\colon \quad \Gamma(W)\to
\Gamma(\bar K W)\,,  
\end{equation}
given by $ D_\omega\psi = (D(\psi e^{\int\omega}))e^{-\int\omega}$ where $\psi
\in\Gamma(W)$, has a non--trivial kernel.  Here $e^{\int
  \omega}\in\Hom(\R^2,\C_*)$ is regarded as a gauge transformation on the
universal cover $\R^2$. Nevertheless, the operator $D_{\omega}$ is still well
defined on the torus $T^2$ because the Leibniz rule of a quaternionic
holomorphic structure implies
\begin{equation}
\label{eq:D_omega}
 D_\omega(\psi)  = D\psi + (\psi \omega)''\,.
\end{equation}
The operator $D_\omega$ is elliptic but due to the term $(\psi \omega)''$ in
(\ref{eq:D_omega}) only a complex linear (rather then quaternionic linear)
operator between the complex rank~2 bundles $W$ and $\bar KW$, where the
complex structure $I$ is given by right multiplication $I(\psi)=\psi i$ by
the quaternion $i$. A section $\psi\in\Gamma(W)$ is in the kernel of
$D_\omega$ if and only if the section $\psi e^{\int \omega}\in\Gamma(\tilde
W)$ is in the kernel of $D$, that is, $\psi e^{\int\omega}\in
H^0(\tilde W)$ is holomorphic. Because the section $\psi e^{\int\omega}$ has
monodromy $h=e^{\int\omega}$, we obtain the $I$--complex linear isomorphism
\begin{equation}
\label{eq:kerD_omega}
\ker D_{\omega}\to H^{0}_{h}(\tilde{W})\colon \quad \psi\mapsto \psi
e^{\int\omega}. 
\end{equation}
In particular, a representation $h=e^{\int\omega}$ is in the spectrum
$\Spec(W,D)$ if and only if $D_{\omega}$ has non--trivial kernel and therefore
\begin{equation}
\label{eq:locus}
\Specs(W,D) = \{ \omega\in\Harm(T^2,\C) \mid \ker D_\omega \neq 0 \}.
\end{equation}

Since $D_{\omega}$ is a holomorphic family of elliptic operators the general
theory of holomorphic families of Fredholm operators (see
Proposition~\ref{prop:fredholm}) implies that the logarithmic spectrum
$\Specs(W,D)$ is a complex analytic subset of $\Hom(\Gamma,\C)\cong \C^2$ and
that $\Spec(W,D)$ is a complex analytic subset of $\Hom(\Gamma,\C_*)\cong
\C_*\times \C_*$. It turns out that the dimension of $\Spec(W,D)\subset
\Hom(\Gamma,\C_*)$ depends on the degree $d=\deg(W)$ of $W$:
\begin{align*}
  \dim(\Spec(W,D))= 2 \qquad \qquad & \textrm{ if } d>0, \\
  \dim(\Spec(W,D))= 1 \qquad \qquad & \textrm{ if } d=0, \\
  \dim(\Spec(W,D))= 0 \qquad \qquad & \textrm{ if } d<0. 
\end{align*}

The case $d\neq 0$ is dealt with in the following lemma.  The case $d=0$ is
the subject of the rest of this paper.

\begin{Lem}
  Let $(W,D)$ be a quaternionic holomorphic line bundle of degree
  $\deg(W)\neq 0$ over a torus~$T^2$. Then
  \begin{enumerate}
  \item[a)] $\Spec(W,D)=\Hom(\Gamma,\C_*)$ if $\deg(W)> 0$ and
  \item[b)] $\Spec(W,D)$ is a finite subset of $\Hom(\Gamma,\C_*)$ if
    $\deg(W)<0$.
  \end{enumerate} 
\end{Lem}

\begin{proof}
  The Fredholm index
  $\Ind(D_\omega)=\dim(\ker(D_\omega))-\dim(\coker(D_\omega))$ (see the proof
  of Lemma~\ref{lem:main1} in Section~\ref{sec:proof_lem_main1} for more
  details) of the elliptic operator $D_\omega$ depends only on the first order
  part. Therefore $\Ind(D_{\omega}$ coincides with the Fredholm index of the
  operator $\dbar$ which, by the Riemann Roch theorem, is given by
  \begin{equation}
    \label{eq:index}
    \Ind(D_\omega) =\Ind(\dbar) = 2 (d-g+1)=2d,
  \end{equation}
  where $d$ denotes the degree of $W$ and the genus $g=1$ for the torus $T^2$.
  The factor $2$ comes from the fact that $\dbar$ is the direct sum of complex
  $\dbar$--operators on $W\cong W_+\oplus W_+$. 
To prove b) we use the quaternionic Pl\"ucker formula with monodromy (see
  appendix to~\cite{BLPP}) which shows that
  \[
  \frac{1}{4\pi} \mathcal{W}(W,D) \geq -d (n+1) + \ord(H)) \] for every
  $(n+1)$--dimensional linear system $H\subset H^0(\tilde W)$ with monodromy.
  This implies that for $d<0$ only a finite number of $h\in \Spec(W,D)$ admit
  non--trivial holomorphic sections with monodromy.
\end{proof}

In case $W=V/L$ is the quaternionic holomorphic line bundle induced by a
conformal immersion $f\colon T^2\rightarrow \S^4$ of a torus into the
conformal 4--sphere the degree $d=\deg(V/L)$ is half of the normal bundle
degree $\deg(\perp_f)$, see \cite{BLPP}.  In particular, for immersions into
the conformal 3--sphere $S^3$ the degree of the induced bundle $V/L$ is always
zero.

\subsection{Spectral curves of degree zero bundles over 2--tori}
Throughout the rest of the paper we assume that $W$ is a quaternionic
holomorphic line bundle of degree zero over a torus $T^2$. The following lemma
shows that in this case the logarithmic spectrum $\Specs(W,D)$, and hence the
spectrum $\Spec(W,D)$, is a $1$--dimensional analytic subset. This allows to
normalize the spectrum to a Riemann surface, the {\em spectral curve}.

\begin{Lem}\label{lem:main1}
  Let $(W,D)$ be a quaternionic holomorphic line bundle of degree zero over a
  torus $T^2$. Then:
\begin{enumerate}[a)]
\item The logarithmic spectrum $\Specs(W,D)$ is a $1$--dimensional analytic
  set in $\Harm(T^2,\C) \cong \C^2$ invariant under translations by
  the lattice $\Gamma^*$. Its normalization is the Riemann surface $\tilde
  \Sigma$ that admits a surjective holomorphic map $\omega\colon \tilde
  \Sigma\rightarrow \Specs(W,D)$ which, on the complement of a discrete set,
  is an injective immersion.
\item The normalization $\tilde \Sigma$ carries a complex holomorphic line
  bundle $\tilde{\mathcal L}$ which is a holomorphic subbundle (in the topology of
  $C^\infty$ convergence) of the trivial $\Gamma(W)$--bundle over $\tilde \Sigma$.
  The fibers of  $\tilde{\mathcal L}$ satisfy  $\tilde{\mathcal
    L}_{\tilde\sigma}\subset\ker D_{\omega(\tilde\sigma)}$ for all
  $\tilde\sigma\in \tilde \Sigma$ with equality away from a discrete set in $\tilde \Sigma$.
\end{enumerate}
\end{Lem}

The fact that $\tilde{\mathcal{L}}$ is a holomorphic line subbundle of the
trivial $\Gamma(W)$--bundle with respect to the Frechet topology of
$C^\infty$--convergence means that a local holomorphic section
$\psi\in H^0(\tilde{\mathcal{L}}_{|U})$ can be viewed as a $C^\infty$--map
\[ 
\psi\colon U\times T^2\rightarrow W\cong T^2\times \C^2 \quad
({\tilde\sigma},p)\mapsto \psi^{\tilde\sigma}(p) 
\] which is holomorphic in ${\tilde\sigma}\in U \subset \tilde\Sigma$.  Here
$W$ is the trivial $\C^2$--bundle with the complex structure given via
quaternionic right multiplication by $i$.

The lemma is proven in Section~\ref{sec:proof_lem_main1} by asymptotic
comparison of the holomorphic families of elliptic operators $D_\omega$ and
$\dbar_\omega$.  This analysis shows that the spectrum $\Spec(W,D)$ of $D$
asymptotically looks like the spectrum $\Spec(\dbar,W)$ of $\dbar$, the {\em
  vacuum spectrum}.

Because the spectrum $\Spec(W,D)=\Specs(W,D)/\Gamma^*$ is the quotient of the
logarithmic spectrum by the dual lattice, Lemma~\ref{lem:main1} implies that
$\Spec(W,D)\subset \Hom(\Gamma,\C_*)$ also is a 1--dimensional analytic
set. The normalization of $\Spec(W,D)$ is the quotient $\Sigma=\tilde
\Sigma/\Gamma^*$ with normalization map $h\colon\Sigma\to\Spec(W,D)$ induced
from $\omega\colon \tilde \Sigma\rightarrow\Specs(W,D)$ via
$h=\exp(\int\omega)$.

\begin{Def}
\label{def:spectral_curve}
Let $(W,D)$ be a quaternionic holomorphic line bundle of degree zero over a
torus. The {\em spectral curve} $\Sigma$ of $W$ is the Riemann surface
normalizing the spectrum $h\colon\Sigma\to\Spec(W,D)$.
\end{Def}
The spectral curve of a conformal immersion $f\colon
T^2\rightarrow S^4$ with trivial normal bundle is \cite{BLPP} the spectral curve of the
induced quaternionic holomorphic line bundle $W=V/L$.  

As proven in Section~\ref{sec:asymptotic_geometry}, the
spectral curve $\Sigma$ of a degree zero quaternionic holomorphic line bundle
$W$ is a Riemann surface with either one or two ends, depending on whether its
genus is infinite or finite. In the former case it is connected, in the latter
case it has at most two components, each containing at least one of the ends.

The bundle $\tilde {\mathcal{L}}$ can be realized as a complex holomorphic
subbundle of the trivial $H^0(\tilde{W})$--bundle (with respect to
$C^{\infty}$--convergence on compact subsets on the universal covering $\R^2$)
over $\tilde \Sigma$ by the embedding
\[ \tilde {\mathcal{L}} \rightarrow H^0(\tilde{W}) \qquad
\psi_{\tilde\sigma}\in \tilde {\mathcal{L}}_{\tilde\sigma} \mapsto
\psi_{\tilde\sigma} e^{\int \omega_{\tilde\sigma}}.\] This line subbundle of
$H^0(\tilde{W})$ is invariant under the action of $\Gamma^*$ and hence
descends to a complex holomorphic line bundle $\mathcal{L}\to \Sigma $. The
fibers of $\mathcal{L}$ satisfy $\mathcal{L}_{\sigma}\subset
H^0_{h^\sigma}(\tilde{W})$ for $\sigma\in\Sigma$ and equality holds away from
a discrete set of points in $\Sigma$.

Recall that the spectrum $\Spec(W,D)\subset \Hom(\Gamma,\C_*)$ is invariant
under the conjugation $\rho(h)=\bar{h}$. The map $\rho$ lifts to an
anti--holomorphic involution $\rho\colon \Sigma\rightarrow \Sigma$ of the
spectral curve $\Sigma$ that satisfies $h\circ \rho=\bar h$.  For every
$\sigma\in \Sigma$ and $\psi\in\mathcal{L}_\sigma$, the section $\psi j$ is a
holomorphic section with monodromy $h^{\rho(\sigma)}=\bar{h^\sigma}$
of~$\tilde W$. For generic $\sigma\in \Sigma$, we have
$\mathcal{L}_{\rho(\sigma)}=H^0_{h^{\rho(\sigma)}}(\tilde W)$ such that $\psi
j\in \mathcal{L}_{\rho(\sigma)}$ and
\[ 
\mathcal{L}_{\rho(\sigma)} = \mathcal{L}_\sigma j\,. 
\] 
By continuity the latter holds for all $\sigma\in \Sigma$.  This shows that
the anti--holomorphic involution $\rho\colon \Sigma\rightarrow \Sigma$ is
covered by right multiplication with $j$ acting on $\mathcal{L}$ and thus has
no fixed points. The following theorem summarizes the discussion so far:

\begin{The}\label{the:main2}
  The spectrum $\Spec(W,D)$ of a quaternionic holomorphic line bundle $(W,D)$
  of degree zero over a torus is a 1--dimensional complex analytic set. Its
  spectral curve, the Riemann surface $\Sigma$ normalizing
  the spectrum $h\colon\Sigma\to\Spec(W,D)$, is equipped with a fixed point
  free anti--holomorphic involution $\rho\colon \Sigma \rightarrow \Sigma$
  satisfying $h\circ \rho = \bar h$.

  There is a complex holomorphic line bundle $\mathcal{L}\to \Sigma$ over
  $\Sigma$, the {\em kernel bundle}, which is a subbundle of the trivial
  bundle $\Sigma\times H^0(\tilde{W})$ with the property that
  $\mathcal{L}_{\sigma} \subset H^{0}_{h^\sigma}(\tilde{W})$ for all
  $\sigma\in \Sigma$ with equality away from a discrete set in $\Sigma$.  The
  bundle $\mathcal{L}$ is compatible with the real structure $\rho\colon
  \Sigma\to\Sigma$ in the sense that
  $\rho^*\mathcal{L}=\mathcal{L}j$. 
  \end{The}

The quotient $\Sigma/\rho$ with respect to the fixed point free involution
$\rho$ is the normalization of the quaternionic spectrum
$\Spec_\H(W,D)$.  If $\Sigma$ is connected, $\Sigma/\rho$ is a
``non--orientable Riemann surface''.

\section{Asymptotic Analysis}\label{sec:analysis}
This section is concerned with a proof of Lemma~\ref{lem:main1} and
some preparatory
results needed in Section~\ref{sec:asymptotic_geometry} to analyze the
asymptotic geometry of the spectrum $\Spec(W,D)$.

The results of Sections~\ref{sec:analysis} and~\ref{sec:asymptotic_geometry}
are obtained by asymptotically comparing the kernels of the holomorphic family
of elliptic operators $D_\omega$ to those of the holomorphic family of
elliptic operators $\dbar_\omega$ arising from the ``vacuum'' $(W,\dbar)$.
For large multipliers, that is, for multipliers in the complement of a compact
subset of $\Hom(\Gamma,\C_*)$, the spectrum $\Spec(W,D)$ is a small
deformation of the vacuum spectrum $\Spec(\dbar,W)$ away from double points of
$\Spec(\dbar,W)$ where handles may form.  The detailed study of the asymptotic
geometry of $\Spec(W,D)$ will then be carried out
Section~\ref{sec:asymptotic_geometry}.

\subsection{Holomorphic families of Fredholm operators and the proof of
  Lemma~\ref {lem:main1}}\label{sec:proof_lem_main1}
The following proposition, combined with
Corollary~\ref{cor:one_dimensionality} proven at the end of this section,
will yield a proof of Lemma~\ref{lem:main1}.

\begin{Pro}\label{prop:fredholm}
  Let $F(\lambda)\colon E_1\rightarrow E_2$ be a holomorphic family of
  Fredholm operators between Banach spaces $E_1$ and $E_2$ parameterized over
  a connected complex manifold $M$. Then the minimal kernel dimension of
  $F(\lambda)$ is attained on the complement of an analytic subset $N\subset
  M$.  If $M$ is 1--dimensional, the holomorphic vector bundle
  $\mathcal{K}_\lambda = \ker(F(\lambda))$ over $M\backslash N$ extends
  through the set~$N$ of isolated points to a holomorphic vector subbundle of
  the trivial $E_1$--bundle over $M$.
\end{Pro}

\begin{proof}
  For $\lambda_0\in M$ there are direct sum decompositions $E_1=E_1' \oplus
  K_1$ and $E_2=E_2' \oplus K_2$ into closed subspaces such that $K_1$ is the
  finite dimensional kernel of $F(\lambda_0)$ and $E_2'$ its image. The
  Hahn--Banach theorem ensures that the finite dimensional kernel of a
  Fredholm operator can be complemented by a closed subspace. The fact that
  the image of an operator of finite codimension $m$ is closed follows from
  the open mapping theorem: for every complement, the restriction of the
  operator to $E_1'$ can be linearly extended to a bijective operator from
  $E_1'\oplus \C^m$ to $E_2=E_2'\oplus K_2$ that respects the direct sum
  decomposition.

  With respect to such a direct sum decomposition, the operators $F(\lambda)$
  can be written as
  \[ F(\lambda) =
  \begin{pmatrix}
    A(\lambda) &  B(\lambda) \\  C(\lambda) &  D(\lambda)  
  \end{pmatrix},
  \] 
  with $A(\lambda)$ invertible for all $\lambda\in U$ contained in an open
  neighborhood $U\subset M$ of $\lambda_0$. For all $\lambda\in U$ we have
  \[  F(\lambda) \underbrace{
  \begin{pmatrix}
    A(\lambda)^{-1} &  -A(\lambda)^{-1} B(\lambda) \\  0  &  \Id_C   
  \end{pmatrix}}_{=:G(\lambda)} = 
  \begin{pmatrix}
    \Id_{E_2'} & 0 \\ C(\lambda) A(\lambda)^{-1} & \underbrace{D(\lambda)
    -C(\lambda)A(\lambda)^{-1} B(\lambda)}_{=:\tilde F(\lambda)}
  \end{pmatrix}.
  \] 
  Because all $G(\lambda)$ are invertible, the kernel of $F(\lambda)$ is
  $G(\lambda)$ applied to $\ker(\tilde F(\lambda))$. Because $\tilde
  F(\lambda)\colon K_1 \rightarrow C$ is an operator between finite dimensional
  spaces whose index coincides with that of $F(\lambda)$, the proposition
  immediately follows from its finite dimensional version: the analytic set
  $N$ is given by the vanishing of all $r\times r$--subdeterminants of $\tilde
  F(\lambda)$, where $r$ is the maximal rank of $\tilde F(\lambda)$ for
  $\lambda\in U$.  For a proof that, in the case of a 1--dimensional parameter
  domain, the bundle of kernels extends through the points in $N$, see e.g.\
  Lemma~23 of~\cite{BFLPP02}.
\end{proof}

\begin{Cor} \label{cor:fredholm}
  If $\Ind(F(\lambda))=0$, the analytic set $N_0=\{\lambda\in M\mid
  \ker(F(\lambda))\neq 0 \}$ can be locally described as the vanishing locus
  of a single holomorphic function.
\end{Cor}
\begin{proof}
  We may assume $N_0\neq M$. Then $N_0=N$, with $N$ as in
  Proposition~\ref{prop:fredholm}, and therefore locally described as the set
  of~$\lambda$ for which $\tilde F(\lambda)$ has a non--trivial kernel.  But
  $\Ind(F(\lambda))=0$ implies that $\tilde F(\lambda)$ is an $r_0\times
  r_0$--matrix with $r_0=\dim(K)=\dim(C)$ such that $N_0$ is locally given as
  the vanishing locus of the determinant of $\tilde F(\lambda)$.
\end{proof}

In order to show that $\Specs(W,D)\subset \Hom(\Gamma,\C)$ is a complex
analytic set, it is sufficient -- by the preceding corollary -- to extend the
family of elliptic operators $D_\omega$ to a holomorphic family of Fredholm
operators defined on a Banach space containing $\Gamma(W)$ such that the
kernels of the extensions coincide with those of~$D_\omega$.
Corollary~\ref{cor:one_dimensionality} below implies that $\Specs(W,D)$ is
neither empty nor the whole space and thus has to be 1--dimensional. Moreover,
by applying Proposition~\ref{prop:fredholm} to $D_\omega$ pulled back by the
normalization $\omega\colon \tilde \Sigma\rightarrow \Specs(W,D)\subset
\Harm(T^2,\C)$, Corollary~\ref{cor:one_dimensionality} shows that the kernel
of $D_\omega$ is 1--dimensional for generic $\omega\in \Specs(W,D)$ and
therefore the kernel bundle is a line bundle over $ \tilde \Sigma$.

\begin{proof}[Proof of Lemma~\ref{lem:main1}.]
  For every integer $k\geq 0$ the operators $D_\omega$ extend to
  bounded operators from the $(k+1)^{th}$--Sobolev space $H^{k+1}(W)$ of
  sections of $W$ to the $k^{th}$--Sobolev space $H^k(\bar KW)$ of sections of
  $\bar KW$. By elliptic regularity (see e.g.\ Theorem~6.30 and Lemma~6.22a of
  \cite{W}), for every $k$ the kernel of the extension to $H^{k+1}(W)$ is
  contained in the space $\Gamma(W)$ of $C^\infty$--sections and therefore
  coincides with the kernel of the original elliptic operator $D_\omega\colon
  \Gamma(W) \rightarrow \Gamma(\bar K W)$.

  The extension of an elliptic operator on a compact manifold to suitable
  Sobolev spaces is always Fredholm: its kernel is finite dimensional (see
  e.g.~\cite{W}, p.258) and so is its cokernel, the space dual to the kernel
  of the adjoint elliptic operator. The Fredholm index
  \[
  \Ind(D_\omega) = \dim(\ker(D_\omega)) - \dim(\coker(D_\omega))
  \]
  of $D_\omega$ depends only on the
  symbol so that by \eqref{eq:index}
  \[
 \Ind(D_\omega) =
  \Ind(\dbar)=2d
  \]
  where  $d$ is the degree of $W$.

  By assumption $d=0$ and therefore Corollary~\ref{cor:fredholm} implies that
  $\Specs(W,D)$ is a complex analytic set locally given as the vanishing locus
  of one holomorphic function in two complex variables.  To see that
  $\Specs(W,D)$ is 1--dimensional it suffices to check that $\Specs(W,D)$ is
  neither empty nor all of $\Harm(T^2,\C)$ which will be proven in
  Corollary~\ref{cor:one_dimensionality}. This corollary also shows that there
  is $\omega\in \Specs(W,D)$ for which $\ker(D_\omega)$ is 1--dimensional.
  Applying Proposition~\ref{prop:fredholm} to $D_\omega$ pulled back to the
  Riemann surface $\tilde \Sigma$ normalizing $\Specs(W,D)$ shows that, for
  generic $\omega\in \tilde\Sigma$, the kernel of $D_\omega$ 1--dimensional.
  The unique extension of $\ker(D_\omega)$ through the isolated points with
  higher dimensional kernel thus defines a holomorphic line bundle
  $\tilde{\mathcal{L}}\rightarrow \tilde \Sigma$.

  Due to elliptic regularity the kernel of the extension of $D_\omega$ to the
  Sobolev space $H^{k+1}(W)$ for every integer $k\geq 0$ is already contained
  in $\Gamma(W)$. In particular, the kernels of the $D_\omega$ do not depend
  on $k$ and define a line subbundle $\tilde{\mathcal{L}}$ of the trivial
  $\Gamma(W)$--bundle over~$\tilde \Sigma$. For every $k\geq 0$ the line
  bundle $\tilde{\mathcal{L}}$ is a holomorphic line subbundle of the trivial
  $H^{k+1}(W)$--bundle and therefore a holomorphic line subbundle of the
  trivial $\Gamma(W)$--bundle where $\Gamma(W)$ is equipped with the Frechet
  topology of $C^\infty$--convergence.
\end{proof}

\subsection{Functional analytic setting} \label{sec:coordinates}
In the following we develop the asymptotic analysis needed to prove
Corollary~\ref{cor:one_dimensionality}. For further applications in
Sections~\ref{sec:asymptotic_geometry} and~\ref{sec:finite_genus} it will be
convenient to view $D_\omega$ as a family of unbounded operators on the Wiener
space of continuous sections of $W$ with absolutely convergent Fourier series
(instead of viewing it as unbounded operators $D_\omega\colon \Gamma(W)\subset
L^2(W)\rightarrow L^2(W)$ which would be sufficient for proving
Corollary~\ref{cor:one_dimensionality}).

In order to work in the familiar framework of function spaces and to simplify
the application of Fourier analysis, we fix a uniformizing coordinate $z$ on
the torus $T^2=\C/\Gamma$ and trivialize the bundle~$W$.  The latter can be
done since the degree of $W$ or, equivalently, of the complex line subbundle
$\hat W=\{v\in W\mid Jv=v i\}$, is zero. To trivialize $W$ we choose a
parallel section $\psi\in \Gamma(\hat W)$ of a suitable flat connection: let
$\hat \nabla$ be the unique flat quaternionic connection on $W$ which is
complex holomorphic, that is, $\hat \nabla J=0$ and $\hat \nabla''=\dbar$, and
which has unitary monodromy $\hat h_\gamma= e^{-\bar b_0 \gamma + b_0 \bar
  \gamma}$ when restricted to $\hat{W}$, where $h\in \Hom(\Gamma,\C_*)$. Then
the connection $\hat \nabla- \bar b_0 dz + b_0 d\bar z$ on $\hat{W}$ is
trivial and we can choose a parallel section $\psi\in \Gamma(\hat W)$ which in
particular satisfies $\dbar \psi = -\psi b_0 dz$.

The choices of uniformizing coordinate $z$ and trivializing section $\psi$
induce isomorphisms
\[ 
C^\infty(T^2,\C^2) \rightarrow \Gamma(W)\,, \qquad  (u_1,u_2) \mapsto
\psi(u_1+j u_2)
\] 
and 
\[ C^\infty(T^2,\C^2)\rightarrow \Gamma(\bar K W)\,, \qquad (u_1,u_2) \mapsto
\psi d\bar z(u_1+j u_2)\,.
\] 
Moreover, via 
\begin{equation}
  \label{eq:coord}
  \omega=(a+\bar b_0) dz+ (b +b_0) d\bar z,
\end{equation} 
 the space $\Harm(T^2,\C)$ is coordinatized by $(a,b)\in \C^2$. Under these
isomorphisms the family of operators $D_\omega$ defined in
\eqref{eq:D_omega_def} takes the form
\begin{equation}
\label{eq:M-eqn}
 D_{a,b}= \bar\partial_{a,b} + M 
 \end{equation}
with
\[
\bar\partial_{a,b}= \begin{pmatrix}
  \frac\partial{\partial\bar z} + b & 0 \\ 
  0 & \frac\partial{\partial z} + a 
\end{pmatrix} 
\qquad \textrm{and} \qquad M=
\begin{pmatrix}
  0 & -\bar q \\ q & 0
\end{pmatrix},
\]
where $q\in C^{\infty}(T^2,\C) $ is the smooth complex function defined by $Q\psi =\psi j dz q$.  We
denote by
\[ \SpeCs=\{ (a,b)\in \C^2 \mid \ker(D_{a,b}) \neq 0\}\] and 
\[\SpeCs_0 =\{
(a,b)\in \C^2 \mid \ker(\bar \partial_{a,b}) \neq 0\}\] the coordinate
versions of the logarithmic spectrum $\Specs(W,D)$ and logarithmic vacuum
spectrum $\Specs(\dbar,W)$.

In the following we equip the space $l^1=l^1(T^2,\C^2)$ of $C^0$--functions
with absolutely convergent Fourier series with the Wiener norm, the
$l^1$--norm
\begin{equation}
  \label{eq:wiener-norm}
  \| u \| = \sum_{c\in \Gamma'} |u_{1,c}| + |u_{2,c}|   
\end{equation}
of the Fourier coefficients of 
\[
u(z)=\left(\sum_{c\in \Gamma'} u_{1,c} \, e^{-\bar c z+c\bar z} , \sum_{c\in
    \Gamma'} u_{2,c} \, e^{-\bar c z+c\bar z}\right) \in l^1(T^2,\C^2), \]
where $\Gamma'\subset \C$ denotes the lattice
\[ \Gamma'=\{c\in \C \mid -\bar c \gamma + c \bar \gamma \in 2\pi i \Z
\textrm{ for all } \gamma\in \Gamma\}.
\] 
The Banach space $l^1$ contains
$C^\infty=C^\infty(T^2,\C^2)$ as a dense subspace (in fact, it contains all
$H^2$--functions, cf.\ the proof of the Sobolev lemma 6.22 in \cite{W}). Since the
inclusion $l^1(T^2,\C^2)\rightarrow C^0(T^2,\C^2)$ is bounded 
$l^1$--convergence implies $C^0$--convergence and in particular pointwise
convergence. We will use those facts in the applications in 
Sections~\ref{sec:asymptotic_geometry} and~\ref{sec:finite_genus}.

\subsection{$D_{a,b}$ and $\dbar_{a,b}$ as closable operators $C^\infty\subset
  l^1\rightarrow l^1$}
We will show  that all operators in the holomorphic families $D_{a,b}$ and
$\dbar_{a,b}$ are closable as unbounded operators $C^\infty\subset
l^1\rightarrow l^1$.  The following
proposition applied to $D_{a,b}$ and
$\dbar_{a,b}$  then allows to compare their kernels by comparing their resolvents.

\begin{Pro}\label{prop:resolvent}
  Let $F\colon D(F)\subset E\rightarrow E$ be a closed operator on a Banach
  space and let $\sigma(F) = \sigma_1 \, \dot\cup\, \sigma_2$ be a
  decomposition of its spectrum into two closed sets such that there exists an
  embedded closed curve $\gamma\colon S^1 \rightarrow \rho(F)=\C\backslash
  \sigma(F)$ enclosing $\sigma_1$. Then the bounded operator
  \[ P:= \frac1{2\pi i} \int_\gamma (\lambda-F)^{-1} d\lambda \]
  is a projection, its range $E_1$ is a subset of $D(F)$, the kernel
  $E_2=\ker(P)$  has $E_2\cap D(F)$ as a dense subspace and both $E_1$
  and $E_2$ are invariant subspaces with $\sigma(F_{|E_i}) = \sigma_i$ for
  $i=1,2$. Moreover, $P$ commutes with every operator that commutes with $F$.
\end{Pro}

For a proof of this proposition see \cite{RS}, Theorems XII.5 and XII.6.

Recall that an unbounded operator $F\colon D(F)\subset E\rightarrow E$ is
closed if the graph norm $\|v\|_F:= \|v\|+\|Fv\|$ is a complete norm on its
domain $D(F)$. For example, if $F$ is a first order elliptic operator $F$, the
extension of the unbounded operator $F\colon C^\infty\subset L^2 \rightarrow
L^2$ to the Sobolev space $H^1$ is closed, because the graph norm is
equivalent to the first Sobolev norm (by Gardings inequality, see e.g.\
\cite{W}, 6.29).  Elliptic regularity (see e.g.\ Theorem~6.30 and Lemma~6.22a
of \cite{W}) shows that this closure $F\colon H^1 \subset L^2 \rightarrow L^2$
has the same kernel as the original operator $F\colon C^\infty\rightarrow
C^\infty$.

In the following we treat $D_{a,b}$ and $\dbar_{a,b}$ as unbounded operators
$C^\infty \subset l^1\rightarrow l^1$ and show (in
Lemma~\ref{lem:D_ab_closable} below) that they also admit closures with kernels
in $C^\infty$.  We will define these closures by taking suitable restrictions
of the respective $L^2$--closures $H^1\subset L^2 \rightarrow L^2$.  For this
we need the following lemma.

\begin{Lem}\label{lem:bounded}
  Let $l^1=l^1(T^2,\C^2)$ be the Banach space of absolutely convergent Fourier
  series with values in~$\C^2$ equipped with the Wiener norm
  \eqref{eq:wiener-norm}. Then:
  \begin{enumerate}
  \item[a)] The multiplication by a $2\times 2$--matrix whose entries are
    functions with absolutely convergent Fourier series (e.g.,
    $C^\infty$--functions) yields a bounded endomorphism of $l^1$ whose
    operator norm is the maximal Wiener norm of the matrix entries.
  \item[b)] A linear operator $\mathcal P\subset l^1\rightarrow l^1$ defined
    on the subspace $\mathcal P$ of Fourier polynomials with the property that
    all Fourier monomials are eigenvectors extends to a bounded operator if and
    only if its eigenvalues are bounded.  Its operator norm is the
    supremum of the moduli of its eigenvalues.
  \end{enumerate}
\end{Lem}

\begin{proof}
  A linear operator $F$ defined on the space of Fourier polynomials with
  values in the Wiener space of absolutely convergent Fourier series uniquely
  extends to a bounded endomorphism of Wiener space if it is bounded on the
  basis $v_k$, $k\in I$, of Fourier monomials of length~1: 
    \[ 
    \|F \sum_k v_k\lambda_k\| \leq \sum_k \| Fv_k\| |\lambda_k| \leq C \sum_k
  |\lambda_k| = C \| \sum_k v_k\lambda_k\|\,,
  \] 
  where  $C$ is such that $\| Fv_k\|\leq C$.
  This proves the claim, because a
  bounded operator defined on a dense subspace of a Banach space has a unique
  extension to the whole space.  The operator norm of the extension is then
  the supremum of the $\| F v_k \|$.
\end{proof}

Let $F$ be the extension to the Sobolev space $H^1$ of one of the operators
$D_{a,b}$ or $\dbar_{a,b}$.  The space $D(F)=\{ v\in l^1 \cap H^1 \mid F(v)\in
l^1 \}\supset C^\infty$ does not depend on the choice of $F$ by Lemma~\ref{lem:bounded}~a)
because the operators $D_{a,b}$
and $\dbar_{a,b}$ differ \eqref{eq:M-eqn} by the bounded endomorphism $M$ of~$l^1$.  
Similarly, when
$F$ is replaced by the extension of another of the operators $D_{a,b}$ and
$\dbar_{a,b}$ in the family, the graph norm of $F\colon D(F)\subset l^1\rightarrow l^1$ is
replaced by an equivalent norm on $D(F)$.  The space $D^1:=D(F)\subset l^1$ is
thus equipped with an equivalence class of norms with respect to which all the
operators $D_{a,b}$ and $\dbar_{a,b}$ extend to bounded operators
$D^1\rightarrow l^1$. In order to see that theses extensions are closed it 
suffices to check that one of the $l^1$--graph norms is complete. We do this
by showing that, for $(a,b)\in \C^2\backslash \SpeCs_0$, the operator
$\dbar_{a,b}$ extends to a bounded operator $\dbar_{a,b}\colon D^1\rightarrow
l^1$ which is injective, surjective, and has a bounded inverse
$(\dbar_{a,b})^{-1}\colon l^1 \rightarrow D^1$.

The kernels of the family of operators $\dbar_{a,b}$---and hence the
logarithmic vacuum spectrum $\SpeCs_0$---are easily understood:
the Schauder basis of $l^1$ given by the Fourier monomials $v_c = ( e^{\bar c
  z - c \bar z} , 0)$ and $w_c = (0, e^{-\bar c z + c \bar z})$ with $c\in
\Gamma'$ is a basis of eigenvectors of all operators $\dbar_{a,b}$, $(a,b)\in
\C^2$, the eigenvalue of $v_c$ being $b-c$ and that of $w_c$ being $a-\bar c$.
In particular, the space of $(a,b)$ for which $\dbar_{a,b}$ has a non--trivial
kernel is
\begin{equation}
  \label{eq:vacuum}
  \SpeCs_0 = (\C\times \Gamma' )\, \cup\, (\bar
  \Gamma'\times \C).
\end{equation}
The kernel of $\dbar_{a,b}$ is 1--dimensional for generic $(a,b)\in \SpeCs_0$.
Exceptions are the double points $(a,b)\in \bar\Gamma' \times \Gamma'$ for
which the kernel is 2--dimensional and the corresponding multiplier
$h_\gamma=e^{(a+\bar b_0)\gamma+(b+b_0)\bar\gamma}$ is real (see also the
discussion at the beginning of Section~\ref{sec:asympt_geo_theorem}).

For $(a,b)\in \C^2\backslash \SpeCs_0$, the operator $\dbar_{a,b}\colon D^1
\rightarrow l^1$ is injective.  Denote by $G_{a,b}$ the unique bounded
endomorphism of $l^1$ extending the operator on the space $\mathcal{P}$ of
Fourier polynomials that is defined by $G_{a,b}(v_c)=(b-c)^{-1} v_c$ and
$G_{a,b}(w_c)=(a-\bar c)^{-1} w_c $.  This unique extension exists by Part~b)
of Lemma~\ref{lem:bounded} (and is compact because it can be approximated by
finite rank operators). The injective endomorphism $G_{a,b}$ maps $l^1$ to the
space $D^1$ and is surjective as an operator $G_{a,b}\colon l^1\rightarrow D^1$,
because $D^1$ is the subspace of $H^1$ of elements
\[u(z)=\sum_{c\in \Gamma'} u_{1,c} v_{-c} + u_{2,c} w_{c}\] for which not only
\[ \sum_{c\in \Gamma'} |u_{1,c}| + |u_{2,c}|<\infty \quad \textrm{ but also }
\quad \sum_{c\in \Gamma'} |b+c||u_{1,c}| + |a-\bar c||u_{2,c}|<\infty.\] In
particular, $G_{a,b}$ is not only a compact endomorphisms of $l^1$, but a
bounded operator $(l^1,\|.\|_{l^1}) \rightarrow (D^1,\|.\|_{\dbar_{a,b}})$,
because the graph norm of $G_{a,b}(u)$ for $u\in l^1$ is
\[
 \|G_{a,b} (u)\|_{\dbar_{a,b}} = \| G_{a,b}(u) \| + \| \dbar_{a,b} G_{a,b}
(u) \| = \| G_{a,b}(u) \| + \|u\|.
\] 
Because
\[
 \dbar_{a,b} G_{a,b}=\Id_{l_1}\,,
 \] 
 the injective operator $\dbar_{a,b}\colon
D^1 \rightarrow l^1$ is also surjective. This proves
\begin{Lem}\label{lem:D_ab_closable}
  For $(a,b)\in \C^2\backslash \SpeCs_0$, the operator $\dbar_{a,b}$ is a
  bijective operator from $D^1$ to $l^1$ with bounded inverse $G_{a,b}=
  \dbar_{a,b}^{-1}\colon l^1\rightarrow D^1$. In particular, for all $(a,b)\in
  \C^2$ the operators $D_{a,b}$ and $\dbar_{a,b}$ extend to closed operators
  $D^1\subset l^1\rightarrow l^1$, because their graph norms on $D^1$ are
  complete.
\end{Lem}

\subsection{Two lemmas about the asymptotics of $D_{a,b}$}
In the asymptotic analysis we make use of an additional symmetry of the
logarithmic vacuum spectrum $\SpeCs_0$, the symmetry induced by the action of
the lattice $\Gamma'$ by quaternionic linear, $J$--commuting gauge
transformations
\[ 
T_c = 
\begin{pmatrix}
  e^{-\bar c z + c \bar z} & 0 \\
  0 & e^{\bar c z - c \bar z} 
\end{pmatrix}, \qquad 
c\in\Gamma' .
\] 
Under this gauge the operators $D_{a,b}=\dbar_{a,b}+M$ and
$\dbar_{a,b}$, $(a,b)\in \C^2$, transform according to
\begin{align}\label{ugly_formula_that_franz_does_not_like}
   D_{a+\bar c,b+c} & = T_c^{-1}(\bar \partial_{a,b} + MT_{-2c}) T_c\,, 
  \\
   \label{not_so_ugly_formula}
   \dbar_{a+\bar c,b+c} & = T_c^{-1}(\bar \partial_{a,b} ) T_c\,. 
\end{align}
The induced $\Gamma'$--action on $\C^2$ is the action of $c\in \Gamma'$ by
\begin{equation}  \label{eq:action_Tc}
  (a,b)\mapsto (a+\bar c,b+c)
\end{equation}
which is a symmetry of the logarithmic vacuum spectrum $\SpeCs_0$, but not of
the logarithmic spectrum$\SpeCs$. In contrast to this, under the $\Gamma'$--action by the gauge
transformations
\[ t_c = \begin{pmatrix}
  e^{-\bar c z + c \bar z} & 0 \\
  0 & e^{-\bar c z + c \bar z}
\end{pmatrix}, \qquad 
c\in\Gamma',
\] 
both $D_{a,b}$ and $\dbar_{a,b}$, $(a,b)\in \C^2$, transform by the same
formula
\begin{align}
   D_{a-\bar c,b+c} & = t_c^{-1}(D_{a,b}) t_c \\
   \dbar_{a-\bar c,b+c} & = t_c^{-1}(\bar \partial_{a,b} ) t_c.
\end{align}
The induced  $\Gamma'$--action on $\C^2$ with
$c\in \Gamma'$ acting by
\begin{equation}
  \label{eq:action_tc}
  (a,b)\mapsto (a-\bar c,b+c)
\end{equation}
is thus a symmetry of both $\SpeCs$ and $\SpeCs_0$. This action is a
coordinate version of the $\Gamma^*$--action in Section~\ref{sec:foundation}
and the corresponding symmetry of $\SpeCs$ and $\SpeCs_0$ is the periodicity
obtained by passing from the spectrum to the logarithmic spectrum.

Instead of $D_{a,b}= \dbar_{a,b} + M$ we consider more generally the
operators $\dbar_{a,b}+MT_{-2c}$ with $c\in \Gamma'$ which also extend to
closed operators $D^1\subset l^1\rightarrow l^1$. If $(a,b)\in \C^2\backslash
\SpeCs_0$, then
\begin{equation}
  \label{eq:factorize}
  \dbar_{a,b}+MT_{-2c} = (\Id + MT_{-2c}G_{a,b}) \dbar_{a,b}
\end{equation}
 where, by Lemma~\ref{lem:bounded}, the operator $\Id + MT_{-2c}G_{a,b}$ is
a bounded endomorphisms of $l^1$. Corollary~\ref{cor:fundamental} below show
that $\Id + MT_{-2c}G_{a,b}$ is invertible if $c$ is sufficiently large.

The following uniform estimate lies at the heart of the asymptotic analysis.

\begin{Lem}\label{lem:asymptotic}
  Let $\Omega\subset \C^2 \backslash \SpeCs_0$ be compact. For every $\delta >0$
  there exists $R>0$ such that
  \[ \| G_{a,b} MT_{-2c} G_{a,b} \| \leq \delta \] for all $(a,b)\in \Omega$ and
  all $c\in \Gamma'$ with $|c|\geq R$.
\end{Lem}
\begin{proof}
  We use that the operators $\bar \partial_{a,b}$ and $G_{a,b}$ are diagonal
  with respect to the Schauder basis $v_c$, $w_c$ of $l^1$.  By Part~b) of
  Lemma~\ref{lem:bounded}, for all $(a,b)\in \Omega$ the norm of the operator
  $G_{a,b}$ satisfies
  \[ \|G_{a,b}\| \leq C:=\frac 1{\min\{C_1,C_2\}}, \] where the minima
  $C_1:=\min\{|a-\bar c|\}>0$ and $C_2:=\min\{|b-c|\}>0$ over all $c\in
  \Gamma'$ and $(a,b)\in \Omega$ are positive, because $\Omega\subset \C^2 \backslash
  \SpeCs_0$ is compact.
 
  Choosing $R'$ such that $|a-\bar c|\geq 3 C \|M\| /\delta$ and $|b- c|\geq 3
  C \|M\| /\delta$ for all $(a,b)\in \Omega$ and $c\in \Gamma'$ with $|c|> R'$
  yields a decomposition $G_{a,b}=G_{a,b}^0+G_{a,b}^{\infty}$ with
  $G_{a,b}^0=G_{a,b}\circ P_{R'}$ and $G_{a,b}^{\infty}= G_{a,b}-G_{a,b}^0$,
  where $P_{R'}$ denotes the projection on $\Span\{ v_c, w_c \mid |c|\leq
  R'\}$ with respect to the splitting $l^1=\Span\{ v_c, w_c \mid |c|\leq
  R'\}\oplus\Span\{ v_c, w_c \mid |c|> R'\}$. Then $\|G_{a,b}^{\infty}\|
  \leq\frac \delta{3C\|M\|}$ and the operator $G_{a,b}^0$ takes values in the
  image of $P_{R'}$ and vanishes on the kernel of $P_{R'}$.

 $M= 
  \begin{pmatrix}
    0 & - \bar q \\ q & 0
  \end{pmatrix}$ with smooth $q$ so that its Fourier series $q(z)=\sum_{c\in \Gamma'} q_c
  e^{-\bar c z + c\bar z}$ converges uniformly and there is $R''$ such that
  $q^\infty = \sum_{c\in \Gamma', |c|>R''} q_c e^{-\bar c z + c\bar z}$ has
  Wiener norm $\sum_{c\in \Gamma', |c|>R''} |q_c| <\frac\delta{3C^2}$. This
  yields the decomposition $M=M^0+M^\infty$ where
  \[ M^\infty = 
  \begin{pmatrix}
    0 & - \bar q_\infty \\ q_\infty & 0
  \end{pmatrix}
  \qquad \textrm{and} \qquad M^0 = 
  \begin{pmatrix}
    0 & - \bar q_0 \\ q_0 & 0
  \end{pmatrix} \] with $q_0 = q-q_\infty$.  By Part~a) of
  Lemma~\ref{lem:bounded}, $\|M^\infty\| \leq \frac \delta{3C^2}$ and the
  operator $M^0$ is the multiplication by a Fourier polynomial.

  For $|c|\geq R:= R'+R''$ we have $G_{a,b}^0 M^0 T_{-2c}G_{a,b}^0=0$ and
  hence, because the shift operator $T_{-2c}$ is an isometry of $l^1$,
  \begin{align*}
    \| G_{a,b} M T_{-2c} G_{a,b} \| & = \| G_{a,b}^\infty M T_{-2c} G_{a,b} +
    G_{a,b}^0 M T_{-2c} G_{a,b}^\infty + G_{a,b}^0 M^\infty T_{-2c} G_{a,b}^0
    \| \leq \delta.
  \end{align*}
\end{proof}

\begin{Cor}\label{cor:fundamental}
  For every $\delta >0$ and every compact $\Omega\subset \C^2\backslash\SpeCs_0$
  there exists $R>0$ such that for all $c\in \Gamma'$ with $|c|\geq R$ and
  $(a,b)\in \Omega$ the operator $\bar
  \partial_{a,b} + MT_{-2c}$ is invertible and \[ \| G_{a,b} - (\bar
  \partial_{a,b} + MT_{-2c})^{-1} \| \leq \delta.\]
\end{Cor}
\begin{proof}
  We use that if an endomorphisms $F$ of a Banach space satisfies $\|F^n\|<1$
  for some power $n$, the operator $\Id +F$ is invertible with inverse given
  by the Neumann series $\sum_{k=0}^\infty (-1)^k F^k$, and
  \begin{equation}
    \label{eq:norm_neumann}
    \|(\Id+F)^{-1}\| \leq \frac1{1-\|F^n\|} (1+ \|F\| + ... + \|F^{n-1}\|).
  \end{equation}
  By Lemma~\ref{lem:asymptotic} we can choose $R>0$ such that, for all for
  $c\in \Gamma'$ with $|c|>R$ and $(a,b)\in \Omega$,
  \[\| G_{a,b} MT_{-2c} G_{a,b} \| \leq
  \min\{\tfrac1{2\|M\|},\tfrac\delta{2(1+G \|M\|)} \},\] where
  $G=\max_{(a,b)\in \Omega} \|G_{a,b}\|$.  Because $T_{-2c}$ is an isometry of
  $l^1$, for $c\in \Gamma'$ with $|c|>R$ and $(a,b)\in \Omega$, the operator
  $MT_{-2c} G_{a,b}$ satisfies $\|(MT_{-2c} G_{a,b})^2\|<1/2$. Thus $\Id
  + MT_{-2c} G_{a,b}$ and therefore $\bar \partial_{a,b} + MT_{-2c}= (\Id +
  MT_{-2c} G_{a,b}) \bar \partial_{a,b}$ are invertible and
  \begin{multline*}
    \| G_{a,b} - (\bar \partial_{a,b} + MT_{-2c})^{-1} \| =
    \| G_{a,b} (\Id - (\Id +   MT_{-2c} G_{a,b})^{-1}) \| \\
    = \| G_{a,b} MT_{-2c} G_{a,b} (\sum_{k=0}^\infty (-1)^k (MT_{-2c}
    G_{a,b})^k) \|\,.
  \end{multline*}
 Together with \eqref{eq:norm_neumann} this implies
  \[ \| G_{a,b} - (\bar \partial_{a,b} + MT_{-2c})^{-1} \| \leq \| G_{a,b}
  MT_{-2c} G_{a,b} \|\frac{1 + \| G_{a,b}\|\|M\|}{1- \|(MT_{-2c}
    G_{a,b})^2\|} <\delta. \]
\end{proof}

The preceding corollary shows in particular that, for every $(a,b)\notin
\SpeCs_0$, the operator $\bar \partial_{a,b} + MT_{-2c}= (\Id + MT_{-2c}
G_{a,b}) \bar \partial_{a,b}$ is invertible if $c\in \Gamma'$ is large enough
and moreover that $(\bar \partial_{a,b} + MT_{-2c})^{-1}$ converges to
$G_{a,b}=(\bar \partial_{a,b})^{-1}$ when $|c|\rightarrow \infty$.  The
convergence is uniform for $(a,b)\in \Omega$ in a compact set $\Omega\subset
\C^2\backslash\SpeCs_0$. This is needed in the proof of the following lemma.

\begin{Lem}\label{lem:projection_tube}
  Let $\delta>0$ and $\Omega\subset \C^2\backslash \SpeCs_0$ compact.  Then there
  exists $R>0$ such that
  \begin{itemize}
  \item[1.)] $\bar\partial_{a,b}+MT_{-2c}$ is invertible for all $(a,b)\in \Omega$
    and all $c\in \Gamma'$ with $|c|>R$.
  \item[2.)] For every ``transversal circle'' $\gamma=\{(a+\lambda,b+\lambda)
    \mid |\lambda|=\epsilon \}$ in $\Omega$ with radius $\epsilon >0$, center $(a,b)\in
    \C^2$, and for all $c\in \Gamma'$ with $|c|>R$, the
    operators
  \[ P^c_\gamma = \frac1{2\pi i} \int_{|\lambda|=\epsilon} (\bar
  \partial_{a+\lambda,b+\lambda} + MT_{-2c})^{-1} d\lambda \qquad \textrm{and}
  \qquad P^\infty_\gamma = \frac1{2\pi i} \int_{|\lambda|=\epsilon}
  G_{a+\lambda,b+\lambda}d\lambda\] are projections and satisfy $\| P^c_\gamma
  - P^\infty_\gamma\| < \delta$.
\end{itemize}
The operator $P^\infty_\gamma$ projects to the finite dimensional sum 
$\im(P^\infty_\gamma)=\bigoplus_{(\tilde a,\tilde b)\in D} \ker(\dbar_{\tilde a,\tilde b})$
 with
$D$ the ``transversal disc'' $D = \{ (a+\lambda,b+\lambda) \in \C^2 \mid
|\lambda|<\epsilon \}$ bounded by the circle~$\gamma$.  If $\delta<1$, for
every $c\in \Gamma'$ with $|c|>R$ the operator $P^c_\gamma$ projects to a
finite dimensional space whose dimension coincides with that of
$\im(P^\infty_\gamma)$ and which is spanned by the kernels of all iterates of
$\dbar_{\tilde a,\tilde b}+MT_{-2c}$ with $(\tilde a,\tilde b)\in D$.
\end{Lem}

The notion \emph{transversal circle} and \emph{transversal disc} reflects the
fact that, away from double points of $\SpeCs_0$, the intersection of a
transversal disc $D$ with $\SpeCs_0$ is transversal.

\begin{proof}
  By Corollary~\ref{cor:fundamental}, there is $R>0$ such that for all
  $(a,b)\in \Omega$ and all $c\in \Gamma'$ with $|c|>R$ the operator $\bar
  \partial_{a,b} + MT_{-2c}$ is invertible and $\|G_{a,b} - (\bar
  \partial_{a,b} + MT_{-2c})^{-1} \| < 2\delta/\diam(\Omega)$.  For every
  transversal circle $\gamma=\{(a+\lambda,b+\lambda) \mid |\lambda|=\epsilon
  \}\subset \Omega$ and every $c\in \Gamma'$ with $|c|>R$, the operators
  $P^c_\gamma$ and $P^\infty_\gamma$ are then well defined and satisfy $\|
  P^c_\gamma - P^\infty_\gamma\| < 2\epsilon \delta/\diam(\Omega) < \delta$.
  Proposition~\ref{prop:resolvent} shows that they are projection operators.

  As one can easily check by evaluation on the Fourier monomials $v_c$ and
  $w_c$, the operator $P^\infty_\gamma$ projects to the space spanned by the
  kernels of $\dbar_{\tilde a,\tilde b}$ for all $(\tilde a,\tilde b)\in D\cap
  \SpeCs_0$.  This space is finite dimensional, because $D\cap \SpeCs_0$ is a
  finite set.

  If $\delta<1$, for every $c\in \Gamma'$ with $|c|>R$ the operator
  $(\Id-P^c_\gamma+P^\infty_\gamma)$ is invertible and maps
  $\im(P^c_\gamma)=\ker(\Id -P^c_\gamma)$ to a subspace of
  $\im(P^\infty_\gamma)$. Similarly, $(\Id-P^\infty_\gamma+P^c_\gamma)$ is
  invertible and maps $\im(P^\infty_\gamma)=\ker(\Id -P^\infty_\gamma)$ to a
  subspace of $\im(P^c_\gamma)$. This shows that $\im(P^c_\gamma)$ and
  $\im(P^\infty_\gamma)$ have the same dimension.  In particular, the space
  $\im(P^c_\gamma)$ is also finite dimensional. By
  Proposition~\ref{prop:resolvent} the finite dimensional image of
  $P^c_\gamma$ is an invariant subspace for the operator
  $\dbar_{a,b}+MT_{-2c}$ and the restriction of $\dbar_{a,b}+MT_{-2c}$ to this
  subset has spectrum contained in $\{|\lambda|<\epsilon\}$. The image of
  $P^c_\gamma$ is thus the direct sum of the kernels of all iterates of
  $\dbar_{\tilde a,\tilde b}+MT_{-2c}$ with $(\tilde a,\tilde b)\in D$.
\end{proof}

\begin{Cor}\label{cor:one_dimensionality}
  There is a point $(a,b)\in \C^2$ such that $\ker(D_{a,b})$ is 0--dimensional
  and there is a point $(\tilde{a},\tilde{b})\in\C^2$ such that
  $\ker(D_{\tilde{a},\tilde{b}})$ is 1--dimensional.
\end{Cor}
\begin{proof}
  Let $(a,b)\in\SpeCs_0$ be a point such that the kernel of
  $\dbar_{a,b}$ is 1--dimensional and choose $\epsilon>0$ such that the
  closure of the transversal disc $D=\{ (a+\lambda,b+\lambda) \in \C^2 \mid
  |\lambda|<\epsilon \}$ intersects $\SpeCs_0$ in $(a,b)$ only. By
  Lemma~\ref{lem:projection_tube} there is $R>0$ such that for all $c\in
  \Gamma'$ with $|c|>R$ there is a unique point $(\tilde a,\tilde b)\in D$ for
  which the kernel of $\dbar_{\tilde a,\tilde b}+MT_{-2c}$ not trivial, but
  1--dimensional.  By \eqref{ugly_formula_that_franz_does_not_like} this
  implies that, for all $c\in \Gamma'$ with $|c|>R$, there is a unique point
  $(\tilde a,\tilde b)\in \SpeCs$ contained in the disc $\{ (a+\bar
  c+\lambda,b+c+\lambda)\in \C^2 \mid |\lambda|<\epsilon\}$ and the kernel of
  $D_{\tilde a,\tilde b}$ at that point $(\tilde a,\tilde b)\in \SpeCs$ is
  1--dimensional.
\end{proof}
Corollary~\ref{cor:one_dimensionality} completes the above proof of
Lemma~\ref{lem:main1} (and hence the proof of Theorem~\ref{the:main2}).

\section{Asymptotic Geometry of Spectral Curves}\label{sec:asymptotic_geometry} 

We investigate the asymptotic geometry of the spectrum $\Spec(W,D)$, the
spectral curve~$\Sigma$, and the kernel bundle $\mathcal{L}\rightarrow \Sigma$
of a quaternionic holomorphic line bundle $(W,D=\dbar+Q)$ of degree zero over a
2--torus.  We show that the spectrum $\Spec(W,D)$ is asymptotic to the vacuum
spectrum $\Spec(W,\dbar)$. As a consequence asymptotically the spectral curve
$\Sigma$ is bi--holomorphic to a pair of planes joined by at most countably
many handles.

\subsection{Statement of the main result}\label{sec:asympt_geo_theorem} 
Recall that the vacuum spectrum
$\Spec(W,\dbar)$ is a real translate of \[\exp(H^0(K))\cup\exp(H^0(\bar K)),\]
see~\eqref{eq:vacuum} or Section~3.2 of \cite{BLPP}. Its double points are the
real representations
\[\Spec(W,\dbar)\cap\Hom(\Gamma,\R_*).\]
If $\exp(\int\omega)=\exp(\int\bar
\eta)\in \Hom(\Gamma,\C_*)$ with $\omega$, $\eta\in H^0(K)$, then
$\omega=\eta$ and $\omega-\bar \omega\in \Gamma^*$ since $1$-forms in the dual lattice are always imaginary.  
On the other hand, for
every $\omega\in \Harm(T^2,\C)$ the representation $\exp(\int\omega)\in
\Hom(\Gamma,\C_*)$ has a unique decomposition $\exp(\int\omega)=
\exp(\int\tfrac{\omega+\bar\omega}2)\exp(\int\tfrac{\omega-\bar\omega}2)$ into
$\R_+$-- and $S^1$--representations the latter of which is real if and only if
$\omega-\bar\omega\in \Gamma^*$, that is, if it is a $\Z_2$--representation.
This shows that the subgroup
\[G'=\{ h \in \Hom(\Gamma,\R_*) \mid h=\exp(\textstyle{\int} \omega)\textrm{
  with } \omega\in H^0(K), \; \omega-\bar \omega \in \Gamma^*\}\] of
$G=\Hom(\Gamma,\C_*)$ acts simply transitive on the set of vacuum double
points.

We now come to the main theorem describing the structure and asymptotics of the spectrum and the kernel line bundle.
\begin{The} \label{the:asymptotic_geo} Let $(W,D)$ be a quaternionic
  holomorphic line bundle of degree zero over a torus.  Then (with respect to
  a fixed bi--invariant metric on $\Hom(\Gamma,\C_*)$):
  \begin{enumerate}[(1)]
  \item For every $\epsilon>0$, there exists a compact set $\Omega\subset
    \Hom(\Gamma,\C_*)$ and a neighborhood $U$ of a vacuum double point in 
    $\Spec(W,\dbar)\cap\Hom(\Gamma,\R_*)$ such that
    \begin{enumerate}[a)]
    \item in the complement of $\Omega$, the spectrum $\Spec(W,D)$ is contained
      in an $\epsilon$--tube around the vacuum spectrum $\Spec(W,\dbar)$;
    \item in the complement of $\Omega$ and away from the neighborhood
      $\bigcup_{h\in G'} h U $ of the vacuum double points, the spectrum
      $\Spec(W,D)$ is a ``graph'' over $\Spec(W,\dbar)$. More precisely,  $\Spec(W,D)$
     is locally a graph over a real translate of $\exp(H^0(K))$ or $\exp(H^0(\bar
      K))$ with respect to coordinates induced, via the exponential map, from
      the splitting $\Hom(\Gamma,\C)\cong \Harm(T^2,\C)=H^0(K)\oplus H^0(\bar
      K)$;
    \item in the neighborhood $h U$, $h\in G'$, of a vacuum double point that
      is contained in the complement of $\Omega$, the intersection of the spectrum
      $\Spec(W,D)$ with $h U$ either consists of one handle, that is, is
      equivalent to an annulus, or it consists of two transversally immersed
      discs which have a double point at a real multiplier.
    \end{enumerate}
    In particular, in the complement of $\Omega$, the spectrum $\Spec(W,D)$ is
    non--singular except at real points which are transversally immersed
    double points. Moreover, for all $h\in \Spec(W,D)\cap
    (\Hom(\Gamma,\C_*)\backslash \Omega)$ either
    \[ \dim(H^0_h(\tilde W))=1 \qquad \textrm{ or } \qquad \dim(H^0_h(\tilde
    W))=2\] depending on whether $h\in (\Hom(\Gamma,\C_*)\backslash
    \Hom(\Gamma,\R))$ or $h\in \Hom(\Gamma,\R)$.
  \item The spectral curve $\Sigma$ is the union
    \[ \Sigma=\Sigma_{cpt}\, \cup \, \Sigma_\infty \] of two $\rho$--invariant
    Riemann surfaces $\Sigma_{cpt}$ and $\Sigma_\infty$ with the following
    properties:
    \begin{enumerate}[a)]
    \item both $\Sigma_{cpt}$ and $\Sigma_\infty$ have a boundary consisting
      of two circles along which they are glued, that is, $\partial
      \Sigma_{cpt} = -\partial \Sigma_{\infty}$;
    \item $\Sigma_{cpt}$ is compact with at most two components each of
      which has a  non--empty boundary; 
    \item $\Sigma_\infty$ consists of two planes, each with a disc removed,
      which are joined by a countable number of handles.
    \end{enumerate}
    In particular, either the spectral curve $\Sigma$ has infinite genus, one
    end and is connected, or it has finite genus, two ends and at most two
    connected components, each containing an end. In the finite genus case,
    both ends are interchanged by the involution $\rho\colon \Sigma
    \rightarrow \Sigma$.
  \item Given $\epsilon>0$ and $\delta>0$, the compact set $\Omega$ and the open
    set $U$ in (1) can be chosen with the following additional properties:
    \begin{enumerate}[a)]
    \item defined on the preimage $V$ under $\tilde \Sigma \rightarrow
      \Spec(W,D)$ of \[\Spec(W,D)\cap \big(\Hom(\Gamma,\C_*)\backslash (\Omega\cup
      \bigcup_{h\in G'} h U)\big),\] there is a holomorphic section $\psi$ of
      $\tilde{\mathcal{L}}\rightarrow \tilde \Sigma$ that is ``$\delta$--close
      to a vacuum solution''. By this we mean that \[ \| \psi - \varphi \|
      <\delta \] for $\varphi$ a nowhere vanishing, locally constant
      section defined over $V$ of the trivial $\Gamma(W)$--bundle over $\tilde
      \Sigma$ that for every $\tilde \sigma\in V\subset \tilde \Sigma$ solves
      $\dbar_\omega\varphi^{\tilde\sigma}=0$ with $\omega\in
      \widetilde\Spec(W,\dbar)$ satisfying $\| \omega(\tilde\sigma)-\omega
      \|<\epsilon$.
    \item Let $h\in G'$ such that the preimage under $\tilde \Sigma
      \rightarrow \Spec(W,D)$ of \[ \Spec(W,D)\cap(\Hom(\Gamma,\C_*)\backslash
      \Omega)\cap h U \] is the sum of two discs. Then $\psi$ holomorphically
      extends through these discs and the extension is again $\delta$--close
      to vacuum solutions.
    \end{enumerate}
  \end{enumerate}
\end{The}

\subsection{Proof of the main result}\label{sec:three_lemma}
Theorem~\ref{the:asymptotic_geo} is a consequence of the following
three lemmas below: Lemma~\ref{lem:tube} shows that, for large multipliers, the
spectrum is contained in an $\epsilon$--tube around the vacuum spectrum with
respect to a bi--invariant metric on $\Hom(\Gamma,\C_*)$.  Lemma~\ref{lem:graph_away_from_dp} shows that, for
large multipliers and away from double points of the vacuum, the spectrum is a
graph over the vacuum spectrum.  Lemma~\ref{lem:double_point} shows that, for
large multipliers, in a neighborhood of a vacuum double point the spectrum
either consists of an annulus or of a pair of discs with a double point.

We continue using the $(a,b)$--coordinates  \eqref{eq:coord} on the Lie algebra
$\Hom(\Gamma,\C)\cong \Harm(T^2,\C)$ of $\Hom(\Gamma,\C_*)$.  In these
coordinates the logarithmic vacuum spectrum is
\begin{equation*}
  \SpeCs_0 = (\C\times \Gamma' )\, \cup\, (\bar
  \Gamma'\times \C)
\end{equation*}
with the set of double points $\bar \Gamma'\times\Gamma'$, see \eqref{eq:vacuum}.
The preimages under the exponential map of points in $\Hom(\Gamma,\C_*)$ are
the orbits of the $\Gamma'$--action $(a,b)\mapsto (a-\bar c,b+c)$ for $c\in \Gamma'$  on $\C^2$, see \eqref{eq:action_tc}. In order to
define fundamental domains for this group action, we fix a basis $c_1$, $c_2$
of $\Gamma'$ of vectors of minimal length. Then both
\begin{align}
  \label{eq:AB}
  \begin{split}
    A & = \{ (a,b) \mid a\in \C \textrm{ and } b=\lambda_1 c_1+ \lambda_2 c_2
    \textrm{ with } \lambda_i\in
    [-1/2,1/2] \} \quad \textrm{ and } \\
    B & = \{ (a,b) \mid b\in \C \textrm{ and } a=\lambda_1 \bar c_1+ \lambda_2
    \bar c_2 \textrm{ with } \lambda_i\in [-1/2,1/2] \}
  \end{split}
\end{align}
are fundamental domains for the $\Gamma'$--action \eqref{eq:action_tc}: orbits
of generic points $(a,b)\in \C^2$ intersect $A$ and $B$ in a single point,
only the orbits of boundary points of $A$ and $B$ intersect several times.
For understanding $\Spec(W,D)\subset \Hom(\Gamma,\C_*)$ it is sufficient to
study the intersection of $\SpeCs$ with the fundamental domain $B$.  To
understand the intersection with $A$ it is sufficient to apply the involution $\rho$, in our coordinates
given by $(a,b)\mapsto(\bar b,\bar a)$, which
interchanges $A$ and $B$ and leaves $\SpeCs$ invariant.

To investigate the intersection $\SpeCs\cap B$ we use 
\eqref{eq:action_Tc} and \eqref{eq:action_tc}: under the $\Gamma'$--action by
the gauge transformation $t_c T_c$ the operator $D_{a,b}$ transforms according
to
\begin{equation}
  \label{eq:gauge_D}
  D_{a,b+2c}= t_c^{-1} T_c^{-1}(\dbar_{a,b}+MT_{-2c}) t_c T_c \qquad
  \textrm{for every} \quad  c\in \Gamma'
\end{equation}
while the fundamental domain $B$ is invariant under the action $(a,b)\mapsto
(a,b+2c)$ of $c\in \Gamma'$.

The following lemma shows that, for large multipliers, the spectrum
$\Spec(W,D)$ lies in an $\epsilon$--tube around the vacuum spectrum
$\Spec(W,\dbar)$.

\begin{Lem}\label{lem:tube}
  For every $\epsilon>0$, there is a compact subset of $B$ in the complement
  of which the intersection of $\SpeCs$ with $B$ is contained in an
  $\epsilon$--tube around $\SpeCs_0$.
\end{Lem}
\begin{proof}
  For $\epsilon >0$ with $2\epsilon< \min\{ |c| \mid c\in
  \Gamma'\backslash\{0\} \} $ we define
  \begin{align}
    \label{eq:Spec_epsilon}
    \SpeCs_0^\epsilon = \{ (a,b) \mid |a-\bar c| < \epsilon \textrm{ for }
    c\in \Gamma' \} \, \cup \, \{ (a,b) \mid |b-c| < \epsilon \textrm{ for }
    c\in \Gamma' \}.
  \end{align}
  For $B'=\{ (a,b)\in B \mid b=\lambda_1 c_1+ \lambda_2 c_2 \textrm{ with }
  \lambda_i\in [-1,1] \}$ we have $B=\bigcup_{c\in \Gamma'} (B'+(0,2c))$.  The
  set $\Omega=B'\backslash \SpeCs_0^\epsilon$ is compact and does not intersect
  $\SpeCs_0$.  Hence, by Corollary~\ref{cor:fundamental}, there is $R>0$ such
  that for all $(a,b)\in \Omega$ and every $c\in \Gamma'$ with $|c|>R$, the kernel
  of $\dbar_{a,b}+MT_{-2c}$ and, by \eqref{eq:gauge_D}, that of $D_{a,b+2c}$
  is zero.  This shows that $\SpeCs$ does not intersect $\bigcup_{c\in
    \Gamma'; |c|> R} (\Omega +(0,2c))$ or, equivalently, that the intersection of
  $\SpeCs$ with the subset $\bigcup_{c\in \Gamma'; |c|> R} (B' +(0,2c))$ of
  $B$ is contained in the $\epsilon$--tube $\SpeCs_0^\epsilon$ around
  $\SpeCs_0$.
\end{proof}

The next lemma and remark show that, away from the double points of the
vacuum spectrum $\Spec(\dbar,W)$, for large enough multipliers the spectrum
$\Spec(W,D)$ is an arbitrarily small deformation of $\Spec(\dbar,W)$.

\begin{Lem} \label{lem:graph_away_from_dp} For $\epsilon>0$ and $\delta>0$
  with $2\epsilon < \min\{ |c| \mid c\in \Gamma'\backslash\{0\}\}$ and
  $\delta<1$ there exists $R>0$ such that:
  \begin{enumerate}[a)]
  \item The intersection of $\SpeCs$ with \[\{|a|<\epsilon\} \times
    \big(\bigcup_{c\in \Gamma'; |c|> R} (\Delta_\epsilon +2c)\big)\subset
    \C^2\] is the graph of a holomorphic function $b\mapsto a(b)$ defined on
    $\bigcup_{c\in \Gamma'; |c|> R} (\Delta_\epsilon +2c)\subset \C$ with
    $\Delta_\epsilon= \{ b=\lambda_1 c_1 + \lambda_2 c_2\mid \lambda_j \in
    [-1,1] \textrm{ and } |b-c|> \epsilon \textrm{ for all } c\in \Gamma'\}$.
    For all points $(a,b)\in \C^2$ contained in this graph, the kernel of
    $D_{a,b}$ is 1--dimensional and, in particular, the resulting multiplier
    is non--real, that is, an element of $\Hom(\Gamma,\C_*)\backslash
    \Hom(\Gamma,\R_*)$.
  \item The bundle $\tilde{\mathcal{L}}\rightarrow \tilde\Sigma$ admits a
    holomorphic section $\psi$ which is defined over the preimage under the
    normalization map $\tilde\Sigma\rightarrow \SpeCs$ of the subset described
    in a) and has the property that, for every $\tilde\sigma$ in this
    preimage, the section $\psi^{\tilde\sigma}\in \tilde{\mathcal L}_{\tilde
      \sigma}$ satisfies
    \[ \| \psi^{\tilde\sigma}-\psi^o\|<\delta,\] where $\|\,\|$ denotes the
    Wiener norm and $\psi^o=(0,1)$ is the Fourier monomial contained in the
    kernel of $\dbar_{0,b}$ for all $b\in \C$.
   \end{enumerate}
\end{Lem}

\begin{Rem}\label{rem:graph_away_from_dp}
  The analogous result for the $a$--plane is obtained by applying the
  anti--holomorphic involution $\rho$: for $\epsilon$, $\delta$, and $R$ as in
  Lemma~\ref{lem:graph_away_from_dp}, the intersection of $\SpeCs$ with the
  set $\bigcup_{c\in \Gamma'; |c|> R} (\bar \Delta_\epsilon +2\bar c) \times
  \{|b|<\epsilon\}$ is the graph of a function $a\mapsto b(a)$ over
  $\bigcup_{c\in \Gamma'; |c|> R} (\bar \Delta_\epsilon +2\bar c)$.  Setting
  $\psi^{\tilde\sigma} = -\psi^{\rho(\tilde\sigma)}j$, the holomorphic
  section $\psi$ from Lemma~\ref{lem:graph_away_from_dp}, b) yields a
  holomorphic section of $\tilde {\mathcal{L}}$ defined on the part of $\tilde
  \Sigma$ which is graph over $\bigcup_{c\in \Gamma'; |c|> R} (\bar
  \Delta_\epsilon +2\bar c)$. This section satisfies
  $\|\psi^{\tilde\sigma}-\psi^\infty\|<\delta$ for $\psi^\infty=-\psi^o j=
  (1,0)$.
\end{Rem}

\begin{proof}
  Let $\tilde \epsilon= \frac15 \epsilon$ and $\Omega=B^\epsilon\backslash
  \SpeCs_0^{\tilde\epsilon}$ with $\SpeCs_0^{\tilde\epsilon}$ as defined in
  \eqref{eq:Spec_epsilon} and
  \[ B^\epsilon=\{ (a,b)\in B \mid \dist((a,b), B') \leq \epsilon\},\] where
  as above $B'=\{ (a,b)\in B \mid b= \lambda_1 c_1+ \lambda_2 c_2 \textrm{
    with } \lambda_i\in [-1,1] \}$. By Lemma~\ref{lem:projection_tube} we can
  chose $R>0$ such that, for every $c\in \Gamma'$ with $|c|>R$, the operator
  $\dbar_{a,b}+MT_{-2c}$ is invertible for all $(a,b)\in \Omega$ and $\|
  P^c_\gamma - P^\infty_\gamma\| < \delta$ for all transversal circles
  $\gamma\subset \Omega$.  As in the proof of Lemma~\ref{lem:tube}, the
  intersection of $\SpeCs$ with $\bigcup_{c\in \Gamma'; |c|> R} (B^\epsilon
  +(0,2c))$ is then contained in the $\tilde \epsilon$--tube
  $\SpeCs_0^{\tilde\epsilon}$ around $\SpeCs_0$.

  For every $b\in \Delta^{4\tilde \epsilon}_{2\tilde \epsilon}$ with
  \begin{equation}
    \label{eq:delta_2epsilon}
     \Delta_{\epsilon_1}^{\epsilon_2} = \{ b\in \C \mid \dist(b,\Delta)\leq
  \epsilon_2 \textrm{ and } |b-c|> \epsilon_1 \textrm{ for all } c\in
  \Gamma'\},
  \end{equation}
  where $\Delta=\{ b=\lambda_1 c_1 + \lambda_2 c_2\mid \lambda_j \in
  [-1,1]\}$, the transversal circle $\gamma_b=\{ (\lambda,b+\lambda) \mid
  |\lambda|=\tilde \epsilon \}$ is contained in $\Omega$.  The operator
  $P^\infty_{\gamma_b}$ projects to the one dimensional kernel of
  $\bar\partial_{0,b}$. Because $\delta<1$,
  Lemma~\ref{lem:projection_tube} implies that, if $c\in \Gamma'$ with
  $|c|>R$, the image of the projection $P^c_\gamma$ is the 1--dimensional
  kernel of $\dbar_{\tilde a,\tilde b}+MT_{-2c}$. Here $(\tilde a,\tilde
  b)\in D$ is the unique point in the transversal disc $D=\{
  (\lambda,b+\lambda) \mid |\lambda|<\tilde \epsilon \}$ for which the kernel
  of $\dbar_{\tilde a,\tilde b}+MT_{-2c}$ is non--trivial.

  By \eqref{eq:gauge_D}, for every $b\in \bigcup_{c\in \Gamma'; |c|> R}
  (\Delta^{4\tilde\epsilon}_{2\tilde\epsilon} +2c)\subset \C$ the disc $\{
  (\lambda,b+\lambda)\in \C^2 \mid |\lambda|<\tilde\epsilon\}$ contains a
  unique point in $\SpeCs$.  This defines a holomorphic function $\lambda$ on
  $\bigcup_{c\in \Gamma'; |c|> R} (\Delta^{4\tilde\epsilon}_{2\tilde\epsilon}
  +2c)$ with $|\lambda|<\tilde \epsilon$ and such that every point in the
  intersection $\SpeCs\cap\big(\{|a|<\epsilon\} \times \bigcup_{c\in \Gamma';
    |c|> R} (\Delta^{3\tilde\epsilon}_{3\tilde\epsilon} +2c)\big)$ is of the form
  $(\lambda(b),b+\lambda(b))$ for some $b\in \bigcup_{c\in \Gamma'; |c|> R}
  (\Delta^{4\tilde\epsilon}_{2\tilde\epsilon} +2c)$.  The Cauchy integral
  formula for the first derivative of $\lambda$ (applied to circles of radius
  $2\tilde \epsilon$) implies that the differential of $\lambda$ restricted to
  $\bigcup_{c\in \Gamma'; |c|> R} (\Delta^{2\tilde\epsilon}_{4\tilde\epsilon}
  +2c)$ is bounded by $1/2$.  Because any two points $b_0$ and $b_1$ in
  $\bigcup_{c\in \Gamma'; |c|> R} (\Delta^{2\tilde\epsilon}_{4\tilde\epsilon}
  +2c)$ can be joined by a curve of length $l\leq \frac\pi2 |b_0-b_1|$, the
  following proposition shows that the map $b\mapsto b+\lambda(b)$ restricted
  to $\bigcup_{c\in \Gamma'; |c|> R}
  (\Delta^{2\tilde\epsilon}_{4\tilde\epsilon} +2c)$ is injective:
\begin{Pro}
  Let $f\colon U\subset \R^n\rightarrow \R^n$ with $\|Df_x-\Id\|\leq \epsilon$
  for $\epsilon>0$ independent of $x\in U$. If any two points $x_0$, $x_1\in
  U$ can be joined by a curve of length $l\leq C|x_0-x_1|$ with $\epsilon
  <1/C$, then $f$ is injective and therefore a diffeomorphism from $U$ to the
  open set $f(U)$.
\end{Pro}
\begin{proof}
  Assume $f(x_0)=f(x_1)$ with $x_0\neq x_1$. Let $\gamma \colon [0,1]
  \rightarrow U$ be a curve of length $l\leq C|x_0-x_1|$ with $\gamma(0)=x_0$,
  $\gamma(1)=x_1$ and constant speed $|\gamma'(t)|=l$ . Then
  \[ |x_0-x_1| = | \int_0^1 (f(\gamma(t))-\gamma(t))' dt | \leq l \int_0^1
  \|Df_{\gamma(t)} -\Id\| dt \leq C|x_0-x_1|\epsilon < |x_0-x_1|.  \] This
  contradicts the assumption that $x_0\neq x_1$ such that $f$ is injective.
  Because $\epsilon<1/C$ and $1/C<1$, the inverse function theorem implies
  that $f$ is a local diffeomorphism.
\end{proof} 
For every $c\in \Gamma'$ with $|c|> R$, the image of the boundary of
$(\Delta^{2\tilde\epsilon}_{4\tilde\epsilon} +2c)$ under the map $b\mapsto
b+\lambda(b)$ is contained in $(\Delta^{3\tilde\epsilon}_{3\tilde\epsilon}
+2c)\backslash (\Delta^{\tilde\epsilon}_{\epsilon} +2c)$ such that the image
of $(\Delta^{2\tilde\epsilon}_{4\tilde\epsilon} +2c)$ is a subset of
$(\Delta^{3\tilde\epsilon}_{3\tilde\epsilon} +2c)$ which contains
$(\Delta^{\tilde\epsilon}_{\epsilon} +2c)$. Therefore, the injective function
$b\mapsto b+\lambda(b)$ maps $\bigcup_{c\in \Gamma'; |c|> R}
(\Delta^{2\tilde\epsilon}_{4\tilde\epsilon} +2c)$ onto a set containing
$\bigcup_{c\in \Gamma'; |c|> R} (\Delta^{\tilde\epsilon}_{\epsilon} +2c)$.
Taking its inverse function $\mu$ restricted to $\bigcup_{c\in \Gamma'; |c|>
  R} (\Delta^{\tilde\epsilon}_{\epsilon} +2c)$ yields a representation of the
intersection of $\SpeCs$ with $\{|a|<\tilde \epsilon\} \times \bigcup_{c\in
  \Gamma'; |c|> R} (\Delta^{\tilde\epsilon}_{\epsilon} +2c)$ as the graph of
the function $a(b)=b-\mu(b)$ over $\bigcup_{c\in \Gamma'; |c|> R}
(\Delta^{\tilde\epsilon}_{\epsilon} +2c)$.

For $b\in\bigcup_{c\in \Gamma'; |c|> R} (\Delta^{\tilde\epsilon}_{\epsilon}
+2c)$ take $b'\in \Delta^{2\tilde\epsilon}_{4\tilde\epsilon}$ and $c\in
\Gamma'$ with $|c|>R$ such that $\mu(b)=b'+2c$.  Then, by definition of
$P^c_{\gamma_{b'}}$ and \eqref{eq:gauge_D},
\begin{equation}
  \label{eq:P_b}
  P_b := t_c^{-1} T_c^{-1}(P^c_{\gamma_{b'}}) t_c T_c =\frac1{2\pi i}
  \int_{|\lambda|=\tilde \epsilon} D_{\lambda,\mu(b)+\lambda}^{-1} d\lambda.
\end{equation}
In particular, the definition of $P_b$ does not depend on the choice of the
representation $\mu(b)=b'+2c$.  Analogously, we define 
\[ 
P_b^\infty = t_c^{-1} T_c^{-1}(P^\infty_{\gamma_{b'}}) t_c T_c=\frac1{2\pi
  i} \int_{|\lambda|=\tilde \epsilon}
{\bar\partial}_{\lambda,\mu(b)+\lambda}^{-1} d\lambda\,.
\] 
This projection
operator is independent of $b$: its kernel contains all Fourier monomials except
the constant section $\psi^o=(0,1)\in C^\infty(T^2,\C^2)$ which spans its
image, that is, $\psi^o=P^\infty_b (\psi^o)$. Since we have chosen $R$ such
that $\|P_b-P_b^\infty\|<\delta$, the section $\psi^{\tilde\sigma}
:=P_b(\psi^o)$ with ${\tilde\sigma}\in \tilde\Sigma$ corresponding to
$(a(b),b)\in \SpeCs$ satisfies $\| \psi^{\tilde\sigma}-\psi^o\|<\delta$.
\end{proof}

The following lemma shows that, for large enough multipliers, in a
neighborhood of a vacuum double point the spectrum $\Spec(W,D)$ either 
has a double point or the vacuum double point resolves into a handle.

\begin{Lem}\label{lem:double_point}
  Let $\epsilon>0$ and $\delta>0$ with $2\epsilon < \min\{ |c| \mid c\in
  \Gamma'\backslash\{0\}\}$ and $\delta<1$. Then there exists $R>0$ such that:
  \begin{enumerate}[a)]
  \item For every $c_0\in \Gamma'$ with $c_0 > 2R$, the intersection of
    $\SpeCs$ with the polydisc \[\{ (a,b)\in \C^2 \mid |a|<\epsilon \textrm{
      and } |b-c_0|<\epsilon \}\] is either bi--holomorphic to an annulus or
    to a pair of transversally intersecting immersed discs with one
    intersection point.  Each disk is a graph over one of the coordinate
    planes.  For a double point $(a,b)$ of the spectrum $\SpeCs$ contained in
    the polydisc, the kernel of the operator $D_{a,b}$ is 2--dimensional and
    the corresponding multiplier is real; for all other $(a,b)\in \SpeCs$
    contained in the polydisc, the kernel of $D_{a,b}$ is 1--dimensional and
    the corresponding multiplier is non--real, that is, an element of
    $\Hom(\Gamma,\C_*)\backslash \Hom(\Gamma,\R_*)$.
  \item The intersection of $\SpeCs$ with $\{|a|<\epsilon\} \times
    \big(\bigcup_{c\in \Gamma'; |c|> R} (\Delta_{\epsilon/2} +2c)\big)\subset
    \C^2$ is the graph of a function $b\mapsto a(b)$ and the holomorphic
    section $\psi$ of $\tilde{\mathcal{L}}$ defined (as in
    Lemma~\ref{lem:graph_away_from_dp}) over this set extends holomorphically
    through the discs around double points which are graphs over the
    $b$--plane. The extension satisfies $\| \psi^{\tilde\sigma}-\psi^o\|<\delta$.
  \end{enumerate}
\end{Lem}

\begin{Rem}\label{rem:double_point}
  As in Remark~\ref{rem:graph_away_from_dp}, the analogous result for the part
  of the spectral curve that is a graph over the $a$--plane can be obtained by
  applying the involution $\rho$.
\end{Rem}

\begin{proof}
  Let $\tilde \epsilon = \frac\epsilon{10}$ and $\Omega=B^{14 \tilde
    \epsilon}\backslash \SpeCs_0^{\tilde\epsilon}$ with $B^\epsilon$ and
  $\SpeCs_0^{\tilde\epsilon}$ as in the proof of
  Lemma~\ref{lem:graph_away_from_dp}.  By Lemma~\ref{lem:projection_tube}, we
  can chose $R>0$ such that for every $c\in \Gamma'$ with $|c|>R$ the operator
  $\dbar_{a,b}+MT_{-2c}$ is invertible for all $(a,b)\in \Omega$ and $\|
  P^c_\gamma - P^\infty_\gamma\| < \delta$ for all transversal circles
  $\gamma\subset \Omega$.  As in the proof of
  Lemma~\ref{lem:graph_away_from_dp}, the intersection of the logarithmic
  spectrum $\SpeCs$ with $\bigcup_{c\in \Gamma'; |c|> R} (B^{14\tilde
    \epsilon} +(0,2c))$ is then contained in the $\tilde \epsilon$--tube
  $\SpeCs_0^{\tilde\epsilon}$ around the vacuum $\SpeCs_0$ and the
  intersection of $\SpeCs$ with $\{|a|<\epsilon\} \times (\bigcup_{c\in
    \Gamma'; |c|> R} (\Delta^{\epsilon}_{\epsilon/2} +2c))$ is a graph over
  $\bigcup_{c\in \Gamma'; |c|> R} (\Delta^{\epsilon}_{\epsilon/2} +2c)$, with
  $\Delta_{\epsilon_1}^{\epsilon_2}$ as in \eqref{eq:delta_2epsilon}.

  For $c_0\in \Gamma'$, not in the same connected component of $\C\backslash
  \bigcup_{c\in \Gamma'; |c|> R} (\Delta^{\epsilon}_{\epsilon/2} +2c)$ than the
  origin, we examine the intersection $\SpeCs\cap\{ (a,b)\in \C^2 \mid
  |a|<\epsilon \textrm{ and } |b-c_0|<\epsilon \}$.  The part of $\SpeCs$
  contained in $\{ (a,b)\in \C^2 \mid |a|<\epsilon \textrm{ and }
  \epsilon/2<|b-c_0|<\epsilon \}$ is then the graph of a function $b\mapsto
  a(b)$ over $\{ b \mid \epsilon/2< |b-c_0|<\epsilon \}$ and, by
  Remark~\ref{rem:graph_away_from_dp} and \eqref{eq:action_tc}, the
  intersection of $\SpeCs$ with $\{ (a,b)\in \C^2 \mid \epsilon/2<|a|<\epsilon
  \textrm{ and } |b-c_0|<\epsilon \}$ is the graph of a function $a\mapsto
  b(a)$ over $\{ a \mid \epsilon/2< |a|<\epsilon \}$.

  We decompose $c_0=c'+2c''$ into $c''\in \Gamma'$ with $|c''|>R$ and $c'
  =l_1c_1+l_2c_2$ for $l_1,\, l_2\in \{0,\pm1\}$.  For $|x|<\frac\epsilon2$
  define the transversal circle
  \begin{equation*}
    \tilde \gamma_x  = \big\{ (0,c')+(x,-x)+ ( \lambda,\lambda ) \mid
    | \lambda| = \frac\epsilon2 + \tilde \epsilon \big \} 
  \end{equation*}
  in $\Omega$ and the corresponding projection operator
  \[ \tilde P_x = t_{c''}^{-1}T_{c''}^{-1}(P^{c''}_{\tilde \gamma_{x}})
  t_{c''} T_{c''}= \frac1{2\pi i} \int_{|\lambda|=\tfrac{\epsilon}2+\tilde
    \epsilon} D_{x+\lambda,c_0-x+\lambda}^{-1} d\lambda. \] Moreover, for $x$
  with $\tilde\epsilon<|x|<\frac\epsilon2$ define the transversal circles
  \begin{align*}
    \gamma^1_x & =\big \{ (0,c')+ (x,-x) + (-x+\mu_1,-x+\mu_1 ) \mid |\mu_1| =
    \tilde \epsilon \big \} \\
    \gamma^2_x & =\big \{ (0,c')+ (x, -x) + (x+\mu_2,x+\mu_2 ) \mid |\mu_2| =
    \tilde \epsilon \big \}
  \end{align*}
  in $\Omega$ and the corresponding projection operators
  \begin{align*}
    P^1_x & = t_{c''}^{-1} T_{c''}^{-1}(P^{c''}_{ \gamma_{x}^1}) t_{c''}
    T_{c''}=\frac1{2\pi i} \int_{|\mu_1|=\tilde \epsilon}
    D_{\mu_1,c_0-2x+\mu_1}^{-1} d\mu_1 \quad \textrm{ and } \\
    P^2_x & = t_{c''}^{-1} T_{c''}^{-1}(P^{c''}_{ \gamma_{x}^2}) t_{c''}
    T_{c''}=\frac1{2\pi i} \int_{|\mu_2|=\tilde \epsilon}
    D_{2x+\mu_2,c_0+\mu_2}^{-1} d\mu_2.
  \end{align*}
  Using the holomorphicity of the resolvent in the definition of $P^c_\gamma$,
  by Stokes theorem we obtain
  \begin{equation}
    \label{eq:sum_of_projectors}
    \tilde P_x = P^1_x + P^2_x
  \end{equation} for all $x$ with $\tilde\epsilon<|x|<\frac\epsilon2$.

  By Lemma~\ref{lem:projection_tube}, for all $|x|<\frac\epsilon2$ the
  operator $\tilde P_x$ projects to a 2--dimensional space which contains the
  span of the kernels of $D_{a,b}$ for all $(a,b) \in \big\{(0,c_0) + (x,-x)+
  (\lambda,\lambda ) \mid |\lambda| < \frac\epsilon2 + \tilde \epsilon \big
  \}$.  For $\tilde\epsilon<|x|<\frac\epsilon2$, the operator $P^1_x$ projects
  to the 1--dimensional kernel of $D_{a,b}$ with $(a,b) \in \big \{ (0,c_0) +
  (x,-x) + (-x+\mu_1,-x+\mu_1 ) \mid |\mu_1| < \tilde \epsilon \big \}$ the
  unique point for which $D_{a,b}$ has a non--trivial kernel.  Analogously,
  $P^2_x$ projects to the 1--dimensional kernel of $D_{a,b}$ for a unique
  $(a,b) \in \big\{ (0,c_0) + (x,-x) + (x+\mu_2,x+\mu_2 ) \mid |\mu_2| <
  \tilde \epsilon \big \}$.

  This gives rise to a holomorphic family of polynomials $p_x(\lambda) =
  \lambda^2 + p_1(x) \lambda + p_2(x)$ (the determinants of the operators
  $D_{x+\lambda,c_0-x+\lambda}$ restricted to the 2--dimensional images of
  $\tilde P_x$) defined on $\{x\mid |x|<\frac\epsilon2\}$ whose zeros describe
  those $\lambda$ with $|\lambda| < \frac\epsilon2 + \tilde \epsilon$ for
  which $\ker(D_{x+\lambda,c_0-x+\lambda})\neq \{0\}$. If
  $\tilde\epsilon<|x|<\frac\epsilon2$, then \eqref{eq:sum_of_projectors}
  implies that, corresponding to the images of $P^1_x$ and $P^2_x$, the
  polynomial $p_x$ has two different zeroes
  \begin{equation}
    \label{eq:two_zeroes}
    \lambda_1(x)=-x+\mu_1(x) \quad \textrm{ and } \quad
    \lambda_2(x)=x+\mu_2(x)     
  \end{equation}
  with $|\mu_k(x)|< \tilde \epsilon$. The discriminant $q(x) = p_1(x)^2 -4
  p_2(x)$ of $p_x$ vanishes exactly at those $x$ for which both zeroes
  coincide.  Its total vanishing order on the set $\{x\mid
  |x|<\frac\epsilon2\}$ is given by the winding number $\frac1{2\pi i}
  \int_{|x|= 2\tilde \epsilon} d(\log(q))$ of $q$ restricted to $|x|= 2\tilde
  \epsilon$.  By \eqref{eq:two_zeroes}, the discriminant
  $q(x)=(\lambda_1(x)+\lambda_2(x))^2-4\lambda_1(x)\lambda_2(x)$ restricted to
  $|x|= 2\tilde \epsilon$ is homotopy equivalent to $4x^2$ (the discriminant
  for the vacuum spectrum). Thus the total vanishing order of the discriminant
  $q$ on the disc $|x|<\epsilon/2$ is two with zeros located in the disc
  $|x|<2\tilde\epsilon$.

  If the discriminant $q$ has two zeros of order one, the intersection of
  $\SpeCs$ with the open set $U=\{(0,c_0) + (x,-x)+ ( \lambda,\lambda ) \mid
  |x|<4\tilde\epsilon, |\lambda| < \frac\epsilon2 + \tilde \epsilon \}$ is
  non--singular and its projection to the disc $\{ |x|<4\tilde\epsilon \}$ is
  a branched 2--fold covering with two branch points of order one. Thus
  $\SpeCs\cap U$ is an annulus whose subsets over
  $3\tilde\epsilon<|x|<4\tilde\epsilon$, by \eqref{eq:two_zeroes}, are
  contained in the sets $\{ (a,b) \mid |a|<\epsilon \textrm{ and }
  \epsilon/2<|b-c_0|<\epsilon \}$ and $\{ (a,b) \mid \epsilon/2<|a|<\epsilon
  \textrm{ and } |b-c_0|<\epsilon \}$ and therefore graphs over the $a$-- and
  $b$--planes. Because $U$ contains the intersection of the polydisc $\{
  (a,b)\in \C^2 \mid |a|\leq\epsilon/2 \textrm{ and } |b-c_0|\leq\epsilon/2
  \}$ with the $\tilde\epsilon$--tube around $\SpeCs_0$, we obtain that the
  intersection of $\SpeCs$ with $\{ (a,b)\in \C^2 \mid |a|<\epsilon \textrm{
    and } |b-c_0|<\epsilon \}$ is an annulus.

  In case the discriminant $q$ has one zero of order two, the intersection of
  $\SpeCs$ with the open set $U=\{(0,c_0) + (x,-x)+ ( \lambda,\lambda ) \mid
  |x|<4\tilde\epsilon, |\lambda| < \frac\epsilon2 + \tilde \epsilon \}$ is a
  2--fold covering of $\{ |x|<4\tilde\epsilon \}$ with one double point over
  the zero of $q$. Thus, the intersection $\SpeCs\cap U$ (and therefore also
  $\SpeCs\cap\{ (a,b)\in \C^2 \mid |a|<\epsilon \textrm{ and }
  |b-c_0|<\epsilon \}$) is normalized by two immersed discs which, near their
  boundaries and hence everywhere, are graphs over the $a$-- and $b$--plane,
  respectively. The two discs intersect transversally, because, by the Cauchy
  integral formula for the first derivative, they are graphs of functions
  $a\mapsto b(a)$ and $b\mapsto a(b)$ with small derivatives.

  I order to see that a double point of $\SpeCs$ in $\{ (a,b)\in \C^2 \mid
  |a|<\epsilon \textrm{ and } |b-c_0|<\epsilon \}$ gives rise to a real
  multiplier, we note that the involution $(a,b)\mapsto (\bar b,\bar a)+(-\bar
  c_0,c_0)$ leaves both $\SpeCs$ and $\{ (a,b)\in \C^2 \mid |a|<\epsilon
  \textrm{ and } |b-c_0|<\epsilon \}$ invariant.  Because there is at most one
  double point of $\SpeCs$ in the polydisc, the double point is a fixed point
  of this involution and hence gives rise to a real multiplier.  The kernel of
  $D_{a,b}$ at the double point is thus 2--dimensional and coincides with the
  image of $\tilde P_x$ for $x$ a zero of the discriminant $q$.

  For all non--singular points $(a,b)\in \SpeCs$ contained in the polydisc $\{
  (a,b)\in \C^2 \mid |a|<\epsilon \textrm{ and } |b-c_0|<\epsilon \}$, the
  kernel of $D_{a,b}$ is 1--dimensional, because the vanishing order of
  $p_x(\lambda)$ seen as a function of two variables is greater or equal to
  the kernel dimension. This completes the proof of part a) of the lemma.

  To prove part b) of the lemma,  we assume that the intersection of
  $\SpeCs$ with the polydisc $\{ (a,b)\in \C^2 \mid |a|<\epsilon \textrm{ and
  } |b-c_0|<\epsilon \}$ consists of two discs with one double point. The
  functions $\lambda_1(x)$ and $\lambda_2(x)$ describing the roots of $p_x$
  for $x\in\{x\mid \tilde\epsilon<|x|<\frac\epsilon2\}$ then extend to
  $\{x\mid |x|<\frac\epsilon2\}$ and define parametrizations of the
  normalization $\tilde \Sigma$ of $\SpeCs$ in a neighborhood of $(0,c_0)$.

  For every $x\in \{x\in \C\mid |x|<\frac\epsilon2\}$, the operator $\tilde
  P_x$ projects to the sum $\tilde{\mathcal{L}}_{\sigma_1(x)} \oplus
  \tilde{\mathcal{L}}_{\sigma_2(x)}$, where $\sigma_j(x)\in \tilde \Sigma$,
  $j=1,2$ are lifts of $(0,c_0) + (x,-x)+ (\lambda_j(x),\lambda_j(x) )\in
  \SpeCs$.  Denote by $\psi_j^{\sigma_j}$, $j=1,2$, nowhere vanishing local
  holomorphic sections of $\tilde {\mathcal{L}}$ defined on $\{ \sigma_j(x)\in
  \tilde \Sigma \mid |x|< \frac\epsilon2\}$.  Then $\tilde P_x\psi^o =
  \psi_1^{\sigma_1(x)} f_1(x) + \psi_2^{\sigma_2(x)} f_2(x)$ for holomorphic
  functions $f_1$, $f_2$, where as in Lemma~\ref{lem:graph_away_from_dp}, we
  set $\psi^o=(0,1)$. By \eqref{eq:sum_of_projectors} we have
  \begin{equation}
    \label{eq:extend}
    P^1_x\psi^o =    \psi_1^{\sigma_1(x)} f_1(x)
  \end{equation}
 for all $x\in \{ x\mid \tilde \epsilon <|x|< \frac\epsilon2\}$. 
  
  Because $x\mapsto c_0-x+\lambda_1(x)$ is bijective near the boundary of $\{
  x\mid |x|< \frac\epsilon2\}$ (which parametrizes a piece of $\SpeCs$ which
  is a graph over the $b$--plane) it is bijective everywhere.  It maps $\{
  x\mid |x|< \frac\epsilon2\}$ onto an open subset of $\{ b\mid |b-c_0|<
  \epsilon + \tilde \epsilon\}$ which contains $\{ b\mid |b-c_0|< 9\tilde
  \epsilon\}$.  Denote by $b\mapsto x(b)$ the inverse of $x\mapsto
  c_0-x+\lambda_1(x)$ restricted to $\{ b\mid |b-c_0|< 9\tilde \epsilon\}$.
  Then $P_b=P^1_{x(b)}$ for all $b\in \{ b\mid 5\tilde \epsilon < |b-c_0|<
  9\tilde \epsilon\}$, where $P_b$ is the operator defined in \eqref{eq:P_b}.
  This allows to extend the holomorphic section $\psi^{\tilde\sigma}=P_b
  \psi^o$ defined in Lemma~\ref{lem:graph_away_from_dp} to the disc
  $\{\tilde\sigma_1(x(b))\in \tilde \Sigma \mid |b-c_0|< 9\tilde \epsilon\}$
  whose image in $\SpeCs$ is a graph over $\{b \mid |b-c_0|< 9\tilde
  \epsilon\}$: for the points ${\tilde\sigma}\in \tilde\Sigma$ over $b$ in the
  annulus $5\tilde \epsilon< |b-c_0|<9\tilde \epsilon$ we have
  $\psi^{\sigma_1(x(b))}=P_b \psi^o=P^1_{x(b)}\psi^o$ such that, by
  \eqref{eq:extend}, the section $\psi_1^{\sigma_1(x(b))} f_1(x(b))$ defines
  an extension to the disc over $\{ b\mid |b-c_0|< 9\tilde \epsilon\}$.  The
  maximum principle implies that this extension still satisfies
  $\|\psi^{\tilde\sigma}-\psi^0\|<\delta$.
\end{proof}

\begin{proof}[Proof of Theorem~\ref{the:asymptotic_geo}.]
  Parts 1) and 3) of the theorem are mere reformulations of
  Lemmas~\ref{lem:tube}, \ref{lem:graph_away_from_dp} and
  \ref{lem:double_point} and Remarks~\ref{rem:graph_away_from_dp} and
  \ref{rem:double_point}.  The decomposition $\Sigma=\Sigma_{cpt}\, \cup \,
  \Sigma_{\infty}$ in Part~2) is also an immediate consequence of
  Lemmas~\ref{lem:graph_away_from_dp} and~\ref{lem:double_point} and
  Remarks~\ref{rem:graph_away_from_dp} and~\ref{rem:double_point}, because the
  spectral curve $\Sigma$ cannot have compact components: on such a compact
  component the harmonic function $\log |h_\gamma|$ had to be constant for all
  $\gamma\in \Gamma$. But this would imply that the normalization map $h\colon
  \Sigma \rightarrow \Hom(\Gamma,\C_*)$ is constant on this component which is
  impossible for an analytic set of dimension one.  Therefore, each component
  of $\Sigma$ contains at least one end for which $h$ goes to infinity and the
  number of components of $\Sigma$ is bounded by the number of ends.

  Asymptotically, away form the vacuum double points, the spectrum
  $\Spec(W,D)$ is a small deformation of the vacuum $\Spec(W,\dbar)$ and
  hence bi--holomorphic to two planes with neighborhoods around the vertices of of $\Z^2$--lattices
  removed. The number of ends of $\Sigma$ depends on the number of handles
  in $\Spec(W,D)$ near large vacuum double points: if there are infinitely many 
  handles the spectral curve $\Sigma$ has one end,  infinite
  genus and is connected. If the number of handles is finite, then $\Sigma$ has
  two ends, at most two components each of which contains an end, and has
  finite genus.
\end{proof}

\subsection{The connection approach to the spectral curve}\label{sec:connection}
The spectral curve $\Sigma$ of a quaternionic holomorphic line bundle $(W,D)$
of degree zero contains a subset $\Sigma_\nabla\subset \Sigma$ that can be
characterized in terms of flat connections adapted to the quaternionic
holomorphic structure $D$. This point of view is advantageous when studying
spectral curves of finite genus.
\begin{Def}
  For a quaternionic holomorphic line bundle $W$ of degree zero over a torus, we 
  define $\Sigma_\nabla\subset \Sigma$ to be the subset of all points
  $\sigma\in\Sigma$ for which non-zero elements $\psi^\sigma\in
  \mathcal{L}_\sigma$ in the fiber over $\sigma$ of the kernel line bundle
  $\mathcal{L}\rightarrow \Sigma$ are nowhere vanishing holomorphic sections
  with monodromy of $W$.
\end{Def}

\begin{Lem}\label{lem:sigma_nabla}
  The subset $\Sigma_\nabla$ is a non--empty open neighborhood of the ends of
  $\Sigma$, that is, the complement $\Sigma\setminus \Sigma_\nabla$ is
  compact.
\end{Lem} 
\begin{proof}
  The fact that $\tilde{\mathcal{L}}$ is a subbundle in the Frechet topology
  of $C^\infty$--convergence implies that every point $\sigma\in
  \Sigma_\nabla$ has a neighborhood in $\Sigma$ on which the non--trivial
  elements of $\mathcal{L}\rightarrow \Sigma$ are nowhere vanishing sections
  with monodromy of $W$. This shows that $\Sigma_\nabla$ is open.

  To see that $\Sigma_\nabla$ is a neighborhood of the ends, note that the
  holomorphic section $\psi$, which has been constructed in
  Lemmas~\ref{lem:graph_away_from_dp} and~\ref{lem:double_point}, is nowhere
  vanishing.  Therefore it is sufficient to check that, for a point $\sigma$
  on a handle joining the two planes in $\Sigma_\infty$ and corresponding to a
  large enough multiplier, a non--trivial section $\psi\in \mathcal{L}_\sigma$
  is nowhere vanishing.  Using the notation in the proof of
  Lemma~\ref{lem:double_point} a non--trivial section $\psi\in
  \tilde{\mathcal{L}}_{\tilde \sigma}$ with $\tilde \sigma\in \tilde \Sigma$
  close to a large vacuum double point $(0,c)$ can be written as $\psi =
  \tilde P_x (\psi^\infty u_a + \psi^o u_b)$ for some point $x$ and $u_a$,
  $u_b\in \C$.  Without loss of generality, we can assume $|u_a|+|u_b|=1$.
  Applying Lemma~\ref{lem:double_point} with $\delta=1/2$ now shows that, in a
  neighborhood of a large enough vacuum double point,
  \[ \| \psi- \psi^\infty u_a - \psi^o u_b\| <\delta (|u_a|+|u_b|)<\tfrac12 \]
  which implies that the section $\psi$ has no zeroes.
\end{proof}

\begin{Rem}
  There are two important special cases of quaternionic holomorphic line
  bundles $(W,D)$ of degree zero for which $\Sigma_\nabla=\Sigma$.  The first
  is the case when the bundle $(W,D)$ carries a flat connection $\nabla$ which
  is a Willmore connection \cite{FLPP01} and adapted to $D$, that is, which
  satisfies the Willmore condition $d^\nabla*Q=0$ and $D=\nabla''$.  The
  spectral curve can then be interpreted as the holonomy eigenline curve of
  the associated family $\nabla^\mu$ of flat connections \cite{FLPP01} which
  means that there is a holomorphic function $\mu\colon \Sigma\rightarrow
  \C_*$, a 2--fold branched covering, such that every non--trivial element
  $\psi^\sigma\in \mathcal{L}_\sigma$ in the fiber of $\mathcal{L}$ over
  $\sigma\in \Sigma$ is a $\nabla^{\mu_\sigma}$--parallel section and hence a
  nowhere vanishing holomorphic section with monodromy of $W$, see
  \cite{FLPP01} or Section~6 of \cite{B}.  Note that trivial Willmore
  connections on a rank~1 bundle correspond to harmonic maps from $T^2$ to
  $S^2$.

  The second is the case when the quaternionic holomorphic line bundle $W$ is
  the bundle $W=V/L$ induced by a conformal immersion $f\colon T^2\rightarrow
  S^4$ with Willmore functional $\mathcal{W}<8\pi$. In this situation
  $\Sigma_\nabla=\Sigma$ is essentially a consequence of Lemma~2.8 of
  \cite{BLPP}.  What remains to be verified is that a non--trivial section
  $\psi^{\sigma_0} \in \mathcal{L}_{\sigma_0}$ over a point ${\sigma_0}\in
  \Sigma$ belonging to the trivial multiplier $h^{\sigma_0}=1$ is nowhere
  vanishing. If such a section $\psi^{\sigma_0}$ had a zero $p$ one could
  construct a 2--dimensional linear system with Jordan monodromy all of whose
  sections vanish at the point $p$ by taking the span of $\psi^{\sigma_0}$ and
  $\frac{\partial \psi^\sigma}{\partial\sigma}_{|\sigma=\sigma_0} + \pi
  \varphi$. Here $\psi^\sigma$ is a local holomorphic section of $\mathcal{L}$
  and $\pi \varphi$ the projection to $W=V/L$ of a parallel section of $V$
  such that $\frac{\partial \psi^\sigma}{\partial\sigma}_{|\sigma=\sigma_0} +
  \pi \varphi$ vanishes at $p$. The quaternionic Pl\"ucker formula with
  monodromy \cite{BLPP} would then contradict $\mathcal{W}<8\pi$.
\end{Rem}

\begin{Def}\label{def:nabla_S}
  For $\sigma\in \Sigma_\nabla$ define the quaternionic connection
  $\nabla^\sigma$ and the complex structure $S^\sigma\in \Gamma(\End(W))$ on~$W$
  by setting
\begin{equation}
  \label{eq:definition_nabla_s}
  \nabla^\sigma\psi^\sigma=0 \quad \textrm{ and } \quad S^\sigma \psi^\sigma
  = \psi^\sigma i\,,
\end{equation} where $\psi^\sigma\in \mathcal{L}_\sigma$ is a non--trivial
element of the fiber $\mathcal{L}_\sigma$ and therefore a nowhere vanishing
holomorphic section with monodromy of $W$.   
\end{Def}
By definition, the connection $\nabla^\sigma$ is flat and compatible with
$S^\sigma$ and $D$, i.e., for $\sigma\in \Sigma_\nabla$
\begin{equation}
  \label{eq:nabla_s_D_compatible}
  \nabla^\sigma
  S^\sigma=0 \qquad \textrm{ and } \qquad D=(\nabla^\sigma)''.  
\end{equation}
The real structure $\rho\colon \Sigma\rightarrow \Sigma$ leaves
$\nabla^\sigma$ invariant and changes the sign of $S^\sigma$, that is,
\begin{equation}
  \label{eq:reality_nabla_s}
  \nabla^{\rho(\sigma)} =\nabla^\sigma \qquad \textrm{ and } \qquad
  S^{\rho(\sigma)} =- S^\sigma.
\end{equation}
By choosing a local holomorphic section $\psi^\sigma$ of the holomorphic line
bundle $\mathcal{L}$ we obtain:
\begin{Lem} \label{lem:nabla_S} The connection $\nabla^\sigma$ and the complex
  structure $S^\sigma$ depend holomorphically on $\sigma\in \Sigma_\nabla$ in
  the sense that
\begin{equation}
  \label{eq:nabla_S_holomorphic}
  (S^\sigma)'= ( S^\sigma)\dot{}\, S^\sigma \qquad \textrm{ and } \qquad
  (\nabla^\sigma)' = (\nabla^\sigma)\dot{}\, S^\sigma,
\end{equation}
where $\dot\,$ and $'$ denote the derivatives with respect to the $t$-- and
$s$--coordinates for $x=t+is$ a local holomorphic chart on $U\subset
\Sigma_\Sigma$.
\end{Lem}

The holomorphic family $S^\sigma$ of complex structures on $W$ defined for
$\sigma \in \Sigma_\nabla$ can be interpreted as a family of holomorphic maps
$S_p\colon \Sigma_\nabla \rightarrow \P_\C(W_p)\cong \CP^1$ parametrized over
the torus $T^2$ .  For this we identify
\[ \P_\C(W_p)\cong\{ S_p\in \End(W_p)\mid S_p^2=-\Id\}\] by identifying the
complex line $v\C$ in $W_p$ with the quaternionic endomorphism $S_p$ whose
$i$--eigenspace is $v\C$.

\begin{The}\label{the:family_of_S}
  Let $(W,D)$ be a quaternionic holomorphic line bundle of degree zero over a
  torus with spectral curve $\Sigma$. For every $p\in T^2$, the evaluation
  $S^\sigma_p$ at $p$ of the complex structure defined in
  \eqref{eq:definition_nabla_s} for all $\sigma\in \Sigma_\nabla$ uniquely
  extends to a holomorphic map
  \begin{equation}
    \label{eq:holomorphic_curve_S}
    S_p\colon \Sigma \rightarrow \P_\C(W_p)\cong \CP^1.
  \end{equation}
  If $\Sigma_\nabla=\Sigma$, this $T^2$--family of holomorphic maps glues to a
  $C^\infty$--map
  \[ S\colon \Sigma\times T^2 \rightarrow \CP^1. \]
\end{The}
\begin{proof}
  Since $\Sigma_\nabla$ is non--empty the evaluation at $p\in T^2$ of a
  non--trivial local holomorphic section $\psi^{\sigma}$ of ${\mathcal{L}}$
  does not vanish identically. Thus we can holomorphically extend
  $\sigma\mapsto \psi^{\sigma}_p\C$ across the isolated zeros of
  $\sigma\mapsto \psi^{\sigma}_p$. This shows that for $p\in T^2$ fixed, $S_p$
  can be holomorphically extended from $\Sigma_\nabla$ to $\Sigma$.

  Recall that $\tilde{\mathcal{L}}$ is a holomorphic subbundle of $\Gamma(W)$
  in the $C^\infty$--topology by Lemma~\ref{lem:main1}. We chose a trivial
  connection on $W$ to identify $\P_\C(W_p)\cong \CP^1$ for all $p\in
  T^2$. Then $\Sigma_\nabla =\Sigma$ implies that $S\colon \Sigma\times T^2\to
  \CP^1$ is smooth.
\end{proof}

If $\Sigma_\nabla\neq \Sigma$, the map $S$ is not necessarily continuous as a
map depending on two variables:  bubbling phenomena
might occur at the points $(\sigma,p)\in \Sigma\times T^2$ for which a
non--trivial $\psi^\sigma\in \mathcal{L}_\sigma$ is a holomorphic section with
monodromy of $W$ with a zero at $p\in T^2$.

\section{Spectral Curves of Finite Genus and the Willmore
  Functional}\label{sec:finite_genus}

We now come to the case where the spectral curve $\Sigma$ of a quaternionic
holomorphic line bundle of degree zero over a torus has finite genus and thus
can be compactified by adding two points $\{o,\infty\}$.
Theorem~\ref{the:convergence} then will show that the $T^2$--family
\eqref{eq:holomorphic_curve_S} of holomorphic maps $S_p\colon \Sigma
\rightarrow \CP^1$ extends to a family of algebraic functions
\[
S_p\colon \Sigma\, \cup \, \{o,\infty\} \rightarrow \CP^1
\]
on the compactification of $\Sigma$. Moreover, the $T^2$--family of complex
holomorphic line bundles corresponding to $S_p$ for $p\in T^2$ move linearly
in the Jacobian of the compactified spectral curve.

Important examples of conformal immersions $f\colon T^2\rightarrow S^4$ with
degree zero normal bundle for which the induced quaternionic holomorphic line
bundle $W=V/L$ has finite spectral genus are constrained Willmore tori with
trivial normal bundle, see~\cite{B}.

\subsection{Asymptotics of finite genus spectral curves}
We say that 
  a quaternionic holomorphic line bundle $W$ of degree zero over a torus has 
   \emph{finite spectral genus} if its spectral curve is of finite
  genus.

In general, the two planes in the
$\Sigma_\infty$--part of the decomposition $\Sigma=\Sigma_{cpt}\, \cup \,
\Sigma_\infty$ of Theorem~\ref{the:asymptotic_geo} are joined by an infinite
number of handles accumulating at the end. In the finite genus case there is a
compact set outside of which there are no handles.  A spectral curve $\Sigma$
of finite genus can thus be compactified by adding two points $\{o,\infty\}$
at infinity.  The real structure $\rho\colon \Sigma\rightarrow \Sigma$ extends
to the compactification $\Sigma\cup\{o,\infty\}$ and interchanges $o$ and
$\infty$.

The compact component $\Sigma_{cpt}$ in the decomposition
$\Sigma=\Sigma_{cpt}\, \cup \, \Sigma_\infty$ of a finite genus spectral curve
can be chosen large enough such that there are no handles joining the two
planes in $\Sigma_\infty$.  Depending on whether $\Sigma_{cpt}$ has one or two
connected components, the compactification $\Sigma \, \cup \, \{o,\infty\}$ is
connected or consists of two connected Riemann surfaces of genus at least one
by
Corollary~\ref{cor:genus_disconnected} below. Since the
$\Sigma_\infty$--component contains no handles, it is the disconnected sum of
two punctured discs,  the neighborhoods of the added points
$\infty$ and $o$, which in the logarithmic picture are  graphs over the $a$-- or
$b$--planes, respectively.  More precisely:

\begin{Lem}\label{lem:parameter_at_infinity}
  Let $W$ be a quaternionic holomorphic line bundle of degree zero over a
  torus $T^2\cong \C/\Gamma$ with spectral curve $\Sigma$ of finite genus.
  Then there is a punctured neighborhood~$U_o$ of one of the points at
  infinity, in the following called $o$,  parameterized by the punctured disc 
  $\{ x\in \C_*\mid |x|<r \}$ for some $r>0$ such that the restriction of the
  normalization map $h\colon \Sigma\rightarrow \Spec(W,D)\subset
  \Hom(\Gamma,\C_*)$ to $U_o$ is of the form
  \[h_\gamma^x = \exp( (\bar b_0 + a(x))\gamma + (b_0 + 1/x)\bar\gamma),\qquad
  \gamma\in \Gamma, \] where $b_0\in \C$ and $x\mapsto a(x)$ is a holomorphic
  function with $a(0)=0$.  Similarly, the other point at infinity, in the
  following called $\infty$, has a punctured neighborhood~$U_\infty$
  parameterized by $\{ x\in \C_*\mid |x|<r \}$ such that the restriction
  of $h$ to $U_\infty$ is
  \[h_\gamma^x = \exp( (\bar b_0 + 1/x)\gamma + (b_0 + b(x))\bar\gamma),\qquad
  \gamma\in \Gamma, \] where $b_0\in \C$ and $x\mapsto b(x)$ is a holomorphic
  function with $b(0)=0$.

  The open sets $U_o$ and $U_\infty$ can be chosen small enough such that the
  sections $\tilde \sigma\mapsto \psi^{\tilde \sigma}$ of
  $\tilde{\mathcal{L}}$ constructed in Lemmas~\ref{lem:graph_away_from_dp}
  and~\ref{lem:double_point} and Remarks~\ref{rem:graph_away_from_dp}
  and~\ref{rem:double_point} are defined on the respective preimages $\tilde
  U_o$, $\tilde U_\infty\subset \tilde \Sigma$ of $U_o$ and $U_\infty$.  By
  setting
  \[ \psi^o=(0,1) \qquad \textrm{ and } \qquad \psi^\infty=(1,0)\] (in the
  trivialization of $W$ used in Sections~\ref{sec:analysis} and
  \ref{sec:asymptotic_geometry}), these sections $\tilde
  \sigma\mapsto\psi^{\tilde \sigma}$ can be extended through the punctures of
  $\tilde U_o$ and $\tilde U_\infty$ such that
  \[ \psi \colon (\tilde U_o \, \cup \, \{o \}) \times T^2 \rightarrow \C^2
  \qquad \textrm{ and } \qquad \psi \colon (\tilde U_\infty \, \cup \,
  \{\infty\}) \times T^2 \rightarrow \C^2
  \]
  are $C^\infty$ as maps depending on two variables and holomorphic in the
  first variable.
\end{Lem}

The main reason for carrying out the asymptotic analysis of
Sections~\ref{sec:analysis} and~\ref{sec:asymptotic_geometry} within the
$l^1$--framework (instead of the usual $L^2$--setting) is that
$l^1$--convergence implies $C^0$--convergence. This is essential for the proof
of Lemma~\ref{lem:parameter_at_infinity}.

\begin{proof} It is sufficient to prove the statement for $U_o$ since the real
  structure $\rho$ exchanges $U_o$ and $U_\infty$ and $\psi^{\rho(\tilde
    \sigma)}:=-\psi^{\tilde \sigma} j$.  In the finite genus case
  Lemmas~\ref{lem:graph_away_from_dp} and~\ref{lem:double_point} imply that,
  for small enough $\delta>0$ and $\epsilon>0$, we can chose $R>0$ big enough
  such that
  \begin{enumerate}[1.)]
  \item the intersection of $\tilde S$ with $\{|a|<\epsilon\} \times
    (\bigcup_{c\in \Gamma'; |c|> R} (\Delta_{\epsilon/2} +2c))$ is a graph of
    a function $b\mapsto a(b)$ over $\bigcup_{c\in \Gamma'; |c|> R}
    (\Delta_{\epsilon/2} +2c)$,
  \item for every $c_0\in \Gamma'$ that satisfies $|c_0|> 2R$, the
    intersection of $\tilde S$ with the polydisc $\{ (a,b)\in \C^2 \mid
    |a|<\epsilon \textrm{ and } |b-c_0|<\epsilon \}$ consists of a pair of
    discs which are graphs over the coordinate planes and have a double point,
    and
  \item the section $\tilde \sigma\mapsto \psi^{\tilde\sigma}$ of $\tilde
    {\mathcal{L}}$ defined by Lemmas~\ref{lem:graph_away_from_dp}
    and~\ref{lem:double_point} over the preimage under the projection $\tilde
    \Sigma\rightarrow \SpeCs$ of the part of $\SpeCs$ which is a graph of a
    function $b\mapsto a(b)$ with $|a|<\epsilon$ over $\bigcup_{c\in \Gamma';
      |c|> R} (B' +2c)$ satisfies
    \begin{equation}
      \label{eq:psi_bounded}
      \|\psi^{\tilde\sigma}-\psi^\infty\| <\delta,
    \end{equation}
    where $B'=\{ (a,b)\in B \mid b=\lambda_1 c_1+ \lambda_2 c_2 \textrm{ with
    } |\lambda_i|\leq 1 \}$.
  \end{enumerate}
  Denote by $\tilde U_o$ an open subset of $\tilde \Sigma$ contained in the
  domain of definition of $\tilde \sigma\mapsto \psi^{\tilde\sigma}$ such that
  the image of $\tilde U_o$ under the projection $\tilde \Sigma\rightarrow
  \SpeCs$ is a graph over $\{ b\in \C\mid |b|>1/r \}$ for some $r>0$. Let
  $U_o$ be the image of $\tilde U_o$ under the projection
  $\tilde\Sigma\rightarrow \Sigma=\tilde \Sigma/\Gamma^*$. By construction,
  this set $U_o$ is a punctured neighborhood of $o$ with the property that the
  restriction of the normalization map $h\colon \Sigma\rightarrow \Spec(W,D)$
  to $U_o$ has a single valued logarithm whose image, in the
  $(a,b)$--coordinates \eqref{eq:coord}, is contained in the $\epsilon$--tube
  around the $b$--plane and is the graph of a holomorphic function $b\mapsto
  a(b)$ which is bounded by $\epsilon$ and defined on $\{ b\in \C\mid |b|>1/r
  \}$.  Setting $x=1/b$ we obtain a parametrization of $U_o$ by $x\in \{x\in
  \C_* \mid |x|<r \}$.  Riemann's removable singularity theorem implies that
  the bounded holomorphic function $x\mapsto a(x)$ extends to $x=0$.  Because
  $\epsilon>0$ can be chosen arbitrarily small, this extension vanishes for
  $x=0$. Using \eqref{eq:coord} this proves that the normalization map~$h$ is
  of the given form when expressed in the $x$--coordinate on~$U_o$ .
  
  Similarly, the section $\tilde \sigma\mapsto \psi^{\tilde\sigma}$ defined on
  $\tilde U_o\subset \tilde \Sigma$, when seen as a holomorphic map from
  $\tilde U_o$ to $C^0(T^2,\C^2)$, is bounded by~\eqref{eq:psi_bounded} and
  Riemann's removable singularity theorem implies that it has a unique
  holomorphic extension to $o$.  This extension maps $o$ to the constant
  element $\psi^o=(0,1)\in C^0(T^2,\C^2)$, because $\delta>0$ can be chosen
  arbitrarily small.  Since by Lemma~\ref{lem:main1} the line bundle
  $\tilde{\mathcal{L}}$ is a holomorphic line subbundle of
  $C^\infty(T^2,\C^2)$ in the $C^\infty$--topology, for every $m\geq 0$ the
  holomorphic section $\tilde \sigma\mapsto \psi^{\tilde\sigma}$ can be seen
  as a holomorphic map from $\tilde U_o$ to $C^m(T^2,\C^2)$ and has a Laurent
  series
  \[ 
  \psi^{\tilde\sigma(x)} = \sum_{k=-\infty}^{\infty} \psi_k^m x^k 
  \] in $C^m(T^2,\C^2)$.  Because in $C^0(T^2,\C^2)$ the section $\tilde
  \sigma\mapsto \psi^{\tilde\sigma}$ has a holomorphic extension through the
  puncture, the Laurent series for $m=0$ is a power series and the
  coefficients of all negative exponents vanish. Since the embedding
  $C^m(T^2,\C^2)\rightarrow C^0(T^2,\C^2)$ is continuous, the uniqueness of
  Laurent series implies that the same is true for all $m$. Thus, for every
  $m\geq 0$ the section $\tilde \sigma\mapsto \psi^{\tilde\sigma}$ extends to
  a holomorphic map from $\tilde U_o\, \cup \,\{o\}$ to $C^m(T^2,\C^2)$ and
  hence $\psi \colon (\tilde U_o \, \cup \, \{o \}) \times T^2 \rightarrow
  \C^2$ is $C^\infty$ and holomorphic in the first variable.
\end{proof}

\begin{Cor} \label{cor:genus_disconnected} Let $\Sigma$ be the spectral curve
 of a quaternionic holomorphic line bundle $(W,D=\dbar+Q)$ of degree
  zero over a torus. Assume $\Sigma$ is disconnected and hence the
  (disconnected) direct sum of two compact Riemann surfaces with a single
  puncture which are interchanged under the anti--holomorphic
  involution~$\rho$.  Then, except in the vacuum case when $Q\equiv 0$, both
  summands have genus $g\geq 1$.
\end{Cor}
It can be shown \cite{B} that the following classes of constrained Willmore
tori in $\S^4$ have irreducible spectral curves: Willmore tori in $S^3$ which
are not M\"obius equivalent to minimal tori in $\R^3$, minimal tori in the
standard 4--sphere or hyperbolic 4--space that are not super--minimal, CMC
tori in $\R^3$ and $\S^3$.

\begin{proof}
  By Theorem~\ref{the:asymptotic_geo} we only have to show that the two
  components have genus $g\geq 1$.  Lemma~\ref{lem:parameter_at_infinity}
  shows that, for each of the components one of the projections which, in the
  $(a,b)$--coordinates of \eqref{eq:coord}, are given by $(a,b)\mapsto a$ and
  $(a,b)\mapsto b$ extends to a non--trivial holomorphic map from the
  compactification of the component onto the torus $\C/\bar\Gamma'$ or
  $\C/\Gamma'$. But by the Riemann--Hurwitz formula a compact surface
  admitting a non--trivial holomorphic map onto a torus has genus $g\geq 1$.
\end{proof}

Using the identification of the Lie algebra $\Hom(\Gamma,\C)$ with
$\Harm(T^2,\C)$ (see Section~\ref{sec:intro}), the logarithmic derivative
$d^\Sigma(\log(h))\in \Omega^1_\Sigma(\Harm(T^2,\C))$ of the normalization map
$h$ can be written as
\begin{equation}\label{eq:def_omegas}
  d^\Sigma(\log(h))= \omega_\infty\, dz+ \omega_o\, d\bar z
\end{equation} 
with $z$ denoting the coordinate induced by the isomorphism $T^2\cong
\C/\Gamma$ used in the definition of the $(a,b)$--coordinates
\eqref{eq:coord}.  The holomorphic forms $\omega_\infty$ and $\omega_o$ are
derivatives $\omega_\infty=da$ and $\omega_o=db$ of the functions $a$ and $b$
which are, up to the $\Gamma'$--action \eqref{eq:action_tc}, well defined on
$\Sigma$. The following corollary is an immediate consequence of
Lemma~\ref{lem:parameter_at_infinity}.

\begin{Cor}
  The form $\omega_\infty$ is holomorphic on $\Sigma\cup\{o\}$ and has a
  second order pole with no residue at $\infty$.  The form
  $\omega_o=\rho^*\bar\omega_\infty$ is holomorphic on $\Sigma\cup\{ \infty\}$
  and has a second order pole with no residue at $o$.
\end{Cor}

The extendibility through the ends of the holomorphic sections $\tilde
\sigma\mapsto \psi^{\tilde\sigma}$ of $\tilde{\mathcal{L}}$ established in
Lemma~\ref{lem:parameter_at_infinity} immediately implies the extendibility of
$S\colon \Sigma_\nabla \rightarrow \Gamma(\End(W))$ through $o$ and $\infty$.
Recall that, on the universal cover of $T^2$, non--trivial elements of
$\tilde{\mathcal{L}}_{\tilde \sigma}$ and $\mathcal{L}_\sigma$ coincide up to
scaling by a complex function if $\tilde \sigma\in \tilde \Sigma$ is the
preimage of $\sigma\in \Sigma$ under the projection $\tilde\Sigma\rightarrow
\Sigma=\tilde \Sigma/\Gamma^*$. Thus $S$ can also be defined using holomorphic
sections of $\tilde{\mathcal{L}}$.

\begin{The}\label{the:convergence}
  Let $W$ be a quaternionic holomorphic line bundle of degree zero over a
  torus with spectral curve $\Sigma$ of finite genus.  By setting
  \[S^{\infty} = J \qquad \textrm{ and } \qquad \S^{o} = -J,\] the family
  \eqref{eq:definition_nabla_s} of complex structures $S^\sigma\in
  \Gamma(\End(W))$ defined for $\sigma\in \Sigma_\nabla$ is extended
  holomorphically (in the $C^\infty$--topology) through the points $o$ and
  $\infty$ to a map
  \[ \sigma\in \Sigma_\nabla\, \cup\, \{o,\infty\} \mapsto S^\sigma\in
  \Gamma(\End(W)).\] In particular, the $T^2$--family
  \eqref{eq:holomorphic_curve_S} of holomorphic functions $S_p \colon
  \Sigma\rightarrow \CP^1$ extends to a family of algebraic functions 
  \[S_p \colon \Sigma\, \cup \, \{o,\infty\} \rightarrow \CP^1, \qquad p\in
  T^2.\] If $\Sigma_\nabla=\Sigma$, this $T^2$--family of algebraic functions
  glues to a $C^\infty$--map
  \[ S\colon (\Sigma\, \cup \, \{o,\infty\})\times T^2\rightarrow \CP^1. \]
\end{The}

In Section~\ref{sec:linear_flow} it will be shown that the complex holomorphic
line bundle belonging to $S_p$, the pull--back of the tautological bundle over
$\CP^1$ by $S_p$, depends linearly on $p\in T^2$ as a map into the Picard
group of the compactified spectral curve.

\subsection{Asymptotics of $\nabla^\sigma$ and the Willmore
  energy}\label{sec:asymptotic_nabla}
We investigate the asymptotics of the connections $\nabla^\sigma$ defined
in \eqref{eq:definition_nabla_s} when $\sigma$ approaches the ends of
$\Sigma$.  Theorem~\ref{the:willmore_energy} shows how the Willmore energy of
a quaternionic holomorphic line bundle of finite spectral genus is encoded in
the asymptotics of $h\colon \Sigma\rightarrow \Hom(\Gamma,\C_*)$.

Because $\rho$ interchanges $o$ and $\infty$ but leaves $\nabla^\sigma$
invariant it is sufficient to investigate $\nabla^\sigma$ in a punctured
neighborhood of $\infty$.  By Theorem~\ref{the:convergence} the sections 
$S^\sigma\in \Gamma(\End(W))$ satisfy $S^\sigma(p)\neq
-J(p)$ for all $p\in T^2$ provided $\sigma$ is in a
small enough punctured neighborhood $U_\infty\subset \Sigma_\nabla$ of
$\infty$. Applying  stereographic projection from $-J$ we write
\begin{equation}
  \label{eq:S_gauge}
  S^\sigma=(1+Y^\sigma) J (1+Y^\sigma)^{-1}
\end{equation}
with $Y^\sigma\in \Gamma(\End_-(W))$. Because $\nabla^\sigma S^\sigma=0$
the flat connection $\nabla^\sigma$ can be expressed as
\begin{equation}
  \label{eq:gauge_nabla}
  \nabla^\sigma= (1+Y^\sigma)\circ (\hat \nabla + \alpha^\sigma) \circ
  (1+Y^\sigma)^{-1}\,.
\end{equation}
Here $\hat\nabla$ denotes the unique flat quaternionic connection with $\dbar
= \hat \nabla''$ and $\hat \nabla J=0$  with unitary holonomy and 
$\alpha^\sigma \in \Omega^1(\End_+(W))$ is a $\hat\nabla$--closed 1--form.  On
the other hand, because $(\nabla^\sigma)''=D=\dbar+Q$, the family
$\nabla^\sigma$ can be written as $\nabla^\sigma = \hat \nabla + Q +
\eta^\sigma$ with $\eta^\sigma\in \Gamma(K\End(W))$ and hence
\begin{align}
  \label{eq:nabla_sigma}
  \nabla^\sigma = \hat \nabla + Q + \eta^\sigma = \hat \nabla +
  \frac1{1+|Y^\sigma|^2} ( \underbrace{\alpha^\sigma + \hat \nabla Y^\sigma
    Y^\sigma - Y^\sigma\alpha^\sigma Y^\sigma}_{\End_+} \underbrace{-\hat
    \nabla Y^\sigma + Y^\sigma\alpha^\sigma - \alpha^\sigma Y^\sigma
  }_{\End_-}).
\end{align}
Since $Q+\eta^\sigma$ has $\End_+$--part of type $K$ we obtain $(\alpha^\sigma
+ \hat \nabla Y^\sigma Y^\sigma - Y^\sigma\alpha^\sigma Y^\sigma)''=0$ and
\begin{multline}\label{eq:QY=alpha_dprime}
  QY^\sigma   = \frac1{1+|Y^\sigma|^2} ( -\hat \nabla Y^\sigma Y^\sigma+
  Y^\sigma\alpha^\sigma Y^\sigma - \alpha^\sigma Y^\sigma Y^\sigma )'' \\
  = \frac1{1+|Y^\sigma|^2} ( \alpha^\sigma - \alpha^\sigma Y^\sigma Y^\sigma
  )'' = (\alpha^\sigma)''.
\end{multline}
The fact that the families $\nabla^\sigma$ and $S^\sigma$ are holomorphic in
$\sigma$ in the sense of \eqref{eq:nabla_S_holomorphic} implies that the
families $Y^\sigma$ and $\alpha^\sigma$ are also holomorphic and satisfy
\[ 
(Y^\sigma)'=(Y^\sigma)\dot{}\, J \qquad \textrm{ and } \qquad
(\alpha^\sigma)'=(\alpha^\sigma)\dot{}\, J. 
\] 
Without loss of generality we
can assume that the punctured neighborhood $U_\infty\subset \Sigma_\nabla$ of
$\infty$ is the parameter domain of a chart $x$ with $x(\infty)=0$ (e.g.\ the
chart defined in Lemma~\ref{lem:parameter_at_infinity}). Then $Y^\sigma$ has a
power series expansion (in the $C^\infty$--topology)
\begin{align}
  \label{eq:asymptotics_W}
  Y^{\sigma(x)} = \sum_{k=1}^\infty Y_k x^k
\end{align}
with $Y_k\in \Gamma(\End_-(W))$ and $(\alpha^\sigma)''$, by
\eqref{eq:QY=alpha_dprime}, has an expansion
\[ 
(\alpha^{\sigma(x)})'' = \sum_{k=1}^\infty \alpha_k'' x^k. 
\]
For every $\sigma\in U_{\infty}$ the form
$\alpha^\sigma\in\Omega^1(\End_+(W))\cong \Omega^1(\C)$ is closed so it has a
unique decomposition $\alpha^\sigma = \alpha^\sigma_{harm} +
\alpha^\sigma_{exact}$ into a harmonic and an exact part.  The multiplier
$h^\sigma$ is then given by $h^\sigma(\gamma)= \hat h(\gamma) e^{-\int_\gamma
  \alpha^\sigma_{harm}}$, where $\hat h$ denotes the holonomy of $\hat \nabla$
restricted to the $i$--eigenline bundle $\hat W$ of $J$ (cf.\
Section~\ref{sec:coordinates}). Lemma~\ref{lem:parameter_at_infinity} now
implies that there is a holomorphic function $a$ on $U_\infty$ with a first
order pole at $\infty$ and a holomorphic function $b$ on
$U_\infty\cup\{\infty\}$ with $b(\infty)=0$ such that
\begin{equation}
  \label{eq:alpha_harm}
  \alpha^\sigma_{harm}=-(a(\sigma)+2\bar b_0) \, dz-b(\sigma)\, d\bar z, 
\end{equation}
where $z$ is the chart on $T^2\cong \C/\Gamma$ and $b_0 \in \C$ satisfies
$\hat h^\gamma = e^{-\bar b_0 \gamma+b_0\bar \gamma}$.

Because both $(\alpha^\sigma)''$ and $(\alpha_{harm}^\sigma)''$ extended
holomorphically through the point $\sigma=\infty$ with $(\alpha^\infty)''=0$
and $(\alpha_{harm}^\infty)''=0$, the same is true for
$(\alpha^\sigma_{exact})''= (\alpha^\sigma-\alpha_{harm}^\sigma)''$.
Moreover, the Fourier expansions of the exact forms $\alpha^\sigma_{exact}$
have no constant terms and the Fourier coefficients of
$(\alpha^\sigma_{exact})'$ and $(\alpha^\sigma_{exact})''$ coincide up to
multiplicative constants independent of $\sigma$ such that
$(\alpha^\sigma_{exact})'$, like $(\alpha^\sigma_{exact})''$, extends
holomorphically through $\infty$ with $(\alpha^\infty_{exact})'=0$.  Since the
two latter components in the decomposition $\alpha^\sigma=
(\alpha^\sigma_{harm})' + (\alpha^\sigma_{exact})' + (\alpha^\sigma)''$ extend
holomorphically through $\sigma=\infty$ and
$(\alpha^\sigma_{harm})'=-(a(\sigma)+2\bar b_0)\, dz$ has a first order pole,
we obtain that $\alpha^\sigma$ has a Laurent series of the form
\begin{align}
  \label{eq:asymptotics_alpha}
  \alpha^{\sigma(x)} = \sum_{k=-1}^\infty \alpha_k x^k 
\end{align}
with closed $\alpha_k \in \Omega^1(T^2,\C)$. The coefficient $\alpha_{-1}$ is
a non--trivial holomorphic 1--form on the torus.  Plugging
\eqref{eq:asymptotics_W} and \eqref{eq:asymptotics_alpha} into
\eqref{eq:nabla_sigma} and taking the $\bar K\End_-(W)$--part yields
$Q=Y_1\alpha_{-1}$. By plugging this into \eqref{eq:QY=alpha_dprime} we obtain
$\alpha_1''=Y_1 \alpha_{-1} Y_1$. Hence
\begin{multline}
  \label{eq:willmore1}
  Res_{x=0} \left(\int_{T^2} \alpha\wedge \frac{\partial\alpha}{\partial
      x}\right)dx = \int_{T^2} \alpha_{-1}\wedge \alpha_1 - \alpha_1 \wedge
  \alpha_{-1} =\\  = -2 \int_{T^2} \alpha_1'' \wedge \alpha_{-1} 
  =  -2 \int_{T^2} Y_1\alpha_{-1} Y_1 \wedge \alpha_{-1} =  i \, \mathcal{W},
\end{multline}
where $\mathcal{W}$ is the Willmore energy of the bundle which is given by 
\[ \mathcal{W} = 2 \int_{T^2} Q\wedge *Q =  2 \int_{T^2}  Y_1\alpha_{-1}\wedge
Y_1*\alpha_{-1}. \]

Using again the identification of Section~\ref{sec:intro} between the Lie
algebra of $\Hom(\Gamma,\C_*)$ and $\Harm(T^2,\C)$, the formula
$h^\sigma(\gamma)= \hat h^\sigma(\gamma) e^{-\int_\gamma
  \alpha^\sigma_{harm}}$ implies $\alpha_{harm}^\sigma =-\log(h^\sigma)+\beta$
for some $\beta\in \Harm(T^2,\C)$ which, like $\log(h)$, is only well
determined up to adding an element of $\Gamma^*$, that is, a $2\pi
i\Z$--periodic harmonic form.  Because $\alpha$ in \eqref{eq:willmore1} can be
replaced by its harmonic part $\alpha_{harm}$, we have proven the following
theorem due to Grinevich and Schmidt, see $(47)$, $(52)$ in \cite{GS98} or
$(44)$ in \cite{Ta06}.

\begin{The}\label{the:willmore_energy}
  Let $(W,D)$ be a quaternionic holomorphic line bundle of degree zero over
  a torus.  In case $(W,D)$ has finite spectral genus its Willmore energy
  is given by
  \[ 
  \mathcal{W} = i Res_{o}\left( \Omega\left( \log(h),d^\Sigma \log(h)
    \right)\right)= - i Res_{\infty}\left(\Omega\left(\log(h),d^\Sigma \log(h)
    \right)\right)\,.
    \] 
 Here $\Omega$ denotes the canonical symplectic form
  \[
  \Omega(\beta_1,\beta_2) := \int_{T^2} \beta_1\wedge \beta_2\,, \qquad
  \beta_1,\beta_2\in \Harm(T^2,\C)\,, 
  \] 
  on the Lie algebra
  $\Hom(\Gamma,\C)\cong \Harm(T^2,\C)$ of $\Hom(\Gamma,\C_*)$ and
  $\log(h)$ denotes the logarithm of $h\colon \Sigma \rightarrow
  \Spec(W,D)\subset \Hom(\Gamma,\C_*)$ which is single valued in punctured
  neighborhoods of $o$ and $\infty$.
\end{The}
Theorem~13.17 in \cite{Hi} is the analogue to
Theorem~\ref{the:willmore_energy} for the energy of harmonic tori
$T^2\rightarrow S^3$ (instead of the Willmore energy of conformal tori
$f\colon T^2\rightarrow S^4$ with finite spectral genus).  To make the analogy
more explicit we give a slight reformulation of the theorem. For this we
define the a skew symmetric product $(\,,\,)_p$ on the space of meromorphic
1--forms with single pole and no residue at $p$ by
\[ (\omega_1,\omega_2)_p = Res_p ( \omega_1F_2), \] where $F_2$ denotes a local
primitive of $\omega_2$, i.e., a holomorphic function with $dF_2=\omega_2$.
Plugging \eqref{eq:def_omegas} into the formula for the Willmore energy we
obtain
\begin{equation}
  \label{eq:residue}
  \mathcal{W} = 4 (\omega_\infty,\omega_o)_\infty \, Vol(\C/\Gamma)= - 4
  (\omega_\infty,\omega_o)_o \, Vol(\C/\Gamma)
\end{equation}
(recall that the forms $\omega_\infty$ and $\omega_o$ defined by
\eqref{eq:def_omegas} depend on the choice of a chart $z$ on $T^2$ which
defines an isomorphism $T^2\cong \C/\Gamma$). For a positive basis $\gamma_1$
and $\gamma_2$ of the lattice $\Gamma$, define $\theta= \omega_\infty \gamma_1
+ \omega_o \bar \gamma_1$ and $\tilde \theta= \omega_\infty \gamma_2 +
\omega_o \bar \gamma_2$.  Because $\gamma_1\bar\gamma_2
-\bar\gamma_1\gamma_2=-2 i\, Vol(\C/\Gamma)$ we obtain
\begin{equation}
  \label{eq:hitchin}
  \mathcal{W} = 2i (\theta,\tilde \theta)_\infty = - 2i (\theta,\tilde
  \theta)_o,
\end{equation}
the direct analogue to the Energy formula given in \cite{Hi}.

As a direct application of \eqref{eq:residue} we show now that the Willmore
energy determines the ``speed'' at which the spectrum $\Spec(W,D)$ converges
to the vacuum $\Spec(W,\dbar)$ when $h$ goes to $\infty$: by
Lemma~\ref{lem:parameter_at_infinity}, a punctured neighborhood of $\infty$
in $\Sigma$ can be parametrized by a parameter $|x|>r$ for which
\[ \log(h^x) = (\bar b_0 + 1/x)\, dz + (b_0 + \lambda x + O(x^2)) \, d\bar
z.\] In this coordinate we thus have $\omega_\infty=da=-1/x^2 dx$ and
$\omega_o=db=(\lambda+O(x)) dx$ such that formula \eqref{eq:residue} implies
\begin{equation}
  \label{eq:convergence_velocity}
  \mathcal{W} = -4 \lambda \, Vol(\C/\Gamma).
\end{equation}

\subsection{The linar flow}\label{sec:linear_flow}
Let $W$ be a quaternionic holomorphic line bundle of degree~$0$ over a 2--torus
$T^2$ with spectral curve $\Sigma$ of finite genus.  The kernel bundle
$\mathcal{L}\rightarrow \Sigma$ does not extend to the compactified spectral
curve $\bar \Sigma=\Sigma\, \cup \, \{ o,\infty\}$ since the monodromy
$h^\sigma$ of elements of $\mathcal{L}_\sigma$ has essential singularities at
$\sigma=o$ and $\infty$.  However, evaluating sections in $\mathcal{L}_\sigma$
at a point $p\in T^2$ gives rise to a complex holomorphic line bundle $E_p\to
\Sigma$, a subbundle of the trivial bundle $\Sigma\times W_p$, which extends
to the compactification $\bar \Sigma$.  Its extension $E_p\to \bar\Sigma$ is
the pull back of the tautological bundle over $\CP^1$ under the algebraic
function $S_p\colon \bar \Sigma\rightarrow \CP^1$ defined in
Theorem~\ref{the:convergence}.  We now prove that the resulting $T^2$--family
of complex holomorphic line bundles $E_p\to \bar\Sigma$ moves linearly in the
Jacobian of $\bar\Sigma$ when the point $p\in T^2$ moves linearly on the
torus.
\begin{The}\label{thm:linearflow}
  Let $W$ be a quaternionic holomorphic line bundle of degree~$0$ over a
  2--torus $T^2$ with spectral curve $\Sigma$ of finite genus and let $p_0\in
  T^2$ be fixed. Then the map
\[
T^2\to \text{Jac}(\bar{\Sigma}):\quad p\mapsto E_p E_{p_0}^{-1}
\]
is a group homomorphism.
\end{The}
\begin{Rem}
  In the special case of a quaternionic holomorphic line bundle that
  corresponds to a harmonic map $f\colon T^2\rightarrow S^2$ from the 2--torus
  to the 2--sphere the above theorem is shown in Chapter~7 of~\cite{Hi}. The
  holomorphic line bundles $E_p$ in that case coincide with the holonomy
  eigenline bundles of the holomorphic family of flat $\SL_2(\C)$--connections
  defined by the harmonic map $f$, cf.\ Section~6.3 of \cite{FLPP01} and
  Section~6.4 of \cite{B}.  To prove the result for general quaternionic
  holomorphic line bundles of degree~$0$ of finite spectral genus (rather than
  bundles corresponding to harmonic maps), we apply similar arguments as in
  Chapter~7 of~\cite{Hi}. In our situation the analog of the harmonic map
  family of flat $\SL_2(\C)$--connections is the family $\nabla^\sigma$ of
  flat quaternionic connections introduced in Section~\ref{sec:connection}.
\end{Rem}
\begin{proof}
Let $V_\infty=U_\infty\; \cup \; \{\infty\}$ be a neighborhood of $\infty$ in
$\bar\Sigma$ with $U_\infty$ as in Section~\ref{sec:asymptotic_nabla} and
denote by $V_o=\rho(V_\infty)$ the corresponding neighborhood of $o$. We
compute the change of $E_p$ in $p\in T^2$ by representing bundles in terms of
\v Cech--cohomology classes with respect to the open cover $\Sigma$,
$V_\infty$ and $V_o$ of $\bar \Sigma$. Denote by $\psi_\infty$ a $\hat
\nabla$--parallel section with monodromy of the quaternionic line bundle
$W\rightarrow T^2$ with complex structure $J$ 
that satisfies $J\psi_\infty=\psi_\infty i$,
where as before $\hat\nabla$ denotes the unique flat connection with
$\hat\nabla''=\dbar$ and $\hat \nabla J=0$ that has unitary holonomy. The
restriction of $E_p$ to $V_\infty$ can then be holomorphically trivialized by
taking the evaluation at $p\in T^2$ of the section
\[ 
\psi^\sigma_\infty := (1+Y^\sigma)\psi_\infty 
\] 
with $Y^\sigma$ as defined
in \eqref{eq:S_gauge}. Similarly, the restriction of $E_p$ to $V_o$ can be
trivialized by taking the evaluation of $\psi^\sigma_o:=
\psi^{\rho(\sigma)}_\infty j$. 

In order to trivialize the restriction of $E_p$ to $\Sigma$, we fix $p_0\in
T^2$ and a nowhere vanishing holomorphic section of the restriction of
$E_{p_0}$ to $\Sigma$. Taking the parallel transport with respect
to~$\nabla^\sigma$, we obtain a family $\psi^\sigma_\Sigma$ of holomorphic
sections of the pullback $\tilde W$ of $W\rightarrow T^2=\C/\Gamma$ to the
universal cover $\C$ of $T^2$ whose restriction $\psi^\sigma_\Sigma(z)$ to
$z\in \C$ is a holomorphic section of $E_p$ for $p\in T^2=\C/\Gamma$, the point
represented by $z$.

The bundle $E_{p_0}$ is represented by the \v Cech--cocycle $f_\infty\colon
U_\infty \rightarrow \C_*$ and $f_o\colon U_o \rightarrow \C_*$ given by $
\psi^\sigma_\Sigma(z_0)= \psi^\sigma_\infty(p_0) f_\infty^\sigma$ for every
$\sigma\in U_\infty = V_\infty\cap \Sigma$ and $\psi^\sigma_\Sigma(z_0) =
\psi^\sigma_o(p_0) f_o^\sigma$ for every $\sigma\in U_o = V_o\cap \Sigma$,
where $z_0\in\C$ denotes a point representing $p_0\in T^2=\C/\Gamma$. Equation
\eqref{eq:gauge_nabla} implies that the bundle $E_p$ is then represented by
the \v Cech--cocycle $\sigma\mapsto f_\infty^\sigma \exp(-\int_{z_0}^{z}
{\alpha}^\sigma)$ and $\sigma\mapsto f_o^\sigma \exp(-\int_{z_0}^{z}
\overline{\alpha^{\rho(\sigma)}})$ defined on $U_\infty$ and $U_o$
respectively for $z\in\C$ a point representing $p\in T^2=\C/\Gamma$. In
Section~\ref{sec:asymptotic_nabla} we have seen that the exact part of
$\alpha^\sigma$ extends holomorphically through $\infty$. Hence $E_{p}$ is
represented by the equivalent cocycle $\sigma\mapsto f_\infty^\sigma
\exp(a(\sigma) (z-z_0))$ and $\sigma\mapsto f_o^\sigma
\exp(\overline{a(\rho(\sigma))(z-z_0)})$ on $U_\infty$ and $U_o$ respectively,
where $a(\sigma)$ is the holomorphic function on $U_\infty$ defined in
\eqref{eq:alpha_harm}.  Changing the point $p\in T^2$ thus amounts to a linear
change of the \v Cech--cohomology class representing $E_p$ which proves the
theorem.
\end{proof}


\end{document}